\documentclass[11pt]{article}
\usepackage[nosumlimits]{amsmath}
\usepackage{amssymb,amsthm,MnSymbol,url}

\newcommand{\pp}[2]{\frac{\partial #1}{\partial #2}}

\newcommand{\del}{\partial}
\def\MM#1{\boldsymbol{#1}}
\DeclareMathOperator{\diff}{d\!}

\usepackage{amscd}
\usepackage{natbib}
\usepackage{helvet}
\usepackage{amsfonts}
\usepackage{algorithmicx}
\usepackage[noend]{algpseudocode}
\usepackage{algorithm}
\newcommand{\jump}[1]{\left[\!\left[ #1 \right]\!\right]}

\usepackage[margin=2cm]{geometry}
\usepackage{graphicx}
\usepackage{color}
\usepackage{caption}
\usepackage{subcaption}
\usepackage{authblk}

\begin{document}
\title{Vertical slice modelling of nonlinear Eady waves using a compatible finite element method}

\author[1,*]{Hiroe Yamazaki}
\author[1]{Jemma Shipton}
\author[2]{Michael J. P. Cullen}
\author[1,3]{Lawrence Mitchell}
\author[1]{Colin J. Cotter}

\affil[1]{Department of Mathematics, Imperial College London, UK}
\affil[2]{Met Office, Exeter, UK}
\affil[3]{Department of Computing, Imperial College London, UK}
\affil[*]{Correspondence to: \texttt{h.yamazaki@imperial.ac.uk}}

\maketitle

\begin{abstract}
A vertical slice model is developed for the Euler-Boussinesq equations with 
a constant temperature gradient in the direction normal to the slice
(the Eady-Boussinesq model). 
The model is a solution of the full three-dimensional equations with no variation 
normal to the slice, which is an idealized problem used to study the formation 
and subsequent evolution of weather fronts. 
A compatible finite element method is used to discretise the governing equations. 
To extend the Charney-Phillips grid staggering in the compatible finite element 
framework, we use the same node locations for buoyancy as the vertical part 
of velocity and apply a transport scheme for a partially continuous finite element 
space. 
For the time discretisation, we solve the semi-implicit equations together with 
an explicit strong-stability-preserving Runge-Kutta scheme to all of the advection 
terms. 
The model reproduces several quasi-periodic lifecycles of fronts despite the 
presence of strong discontinuities. 
An asymptotic limit analysis based on the semi-geostrophic theory shows that 
the model solutions are converging to a solution in cross-front 
geostrophic balance. 
The results are consistent with the previous results using finite 
difference methods, indicating that the compatible finite element method is 
performing as well as finite difference methods for this test problem. 
We observe dissipation of kinetic energy of the cross-front velocity 
in the model due to the lack of resolution at the fronts, 
even though the energy loss is not likely to account for the large gap on the strength 
of the fronts between the model result and the semi-geostrophic limit solution.
\end{abstract}

\noindent \textbf{keywords:} mixed finite elements; frontogenesis; Eady model; 
asymptotic convergence; semi-geostrophic; numerical weather prediction

\section{Introduction}

In the last two decades, Finite element methods have become a popular
discretisation approach for numerical weather prediction (NWP). 
The main focus has been on spectral elements or discontinuous Galerkin (DG) 
methods \citep{fournier2004spectral, thomas2005ncar, dennis2012cam, 
kelly2012continuous, kelly2013implicit, marras2013simulations, 
brdar2013comparison, bao2015horizontally, marras2015review}. 
Another track of research, which we continue here, has been on compatible
finite element methods \citep[e.g.][]{cotter2012mixed, staniforth2013analysis, 
cotter2014finite, mcrae2014energy, natale2016compatible}.
This work is motivated by the need to move away from the conventional
latitude-longitude grids, whilst
retaining properties of the Arakawa C-grid staggered finite difference method.
The compatible finite element method is a family of mixed finite element 
methods where different finite element 
spaces are selected for different variables. 
Compatible finite element methods are built from finite element spaces that have 
differential operators such as grad and curl that map from one space to another. 
This embedded property makes the compatible finite element method analogous 
to the Arakawa C-grid staggered finite difference method 
\citep[][]{arakawa1977computational} with extra flexibility in the choice of 
discretisation to optimize the ratio between global velocity degrees of 
freedom (DoFs) and global pressure DoFs for the sake of avoiding spurious 
modes \citep[][]{cotter2012mixed}.
In addition, it allows the use of arbitrary grids with no requirement of 
orthogonality.

In this study, the compatible finite element method is applied to the 
Euler-Boussinesq equations with a constant temperature gradient in the 
$y$-direction: the Eady-Boussinesq vertical slice model 
\citep[][]{hoskins1972atmospheric}. 
The model is a solution of the full three-dimensional equations with no 
variation normal to the slice. 
As the domain of the vertical slice model consists of 
a two-dimensional slice, 
it can be run much quicker on a workstation than a full three-dimensional model.
This makes the model ideal for numerical studies and test problems for 
NWP models.
There have been many studies on this idealized problem to examine the 
formation and subsequent evolution of weather fronts 
\citep[e.g.][]{williams1967atmospheric, nakamura1989nonlinear, 
nakamura1994nonlinear, snyder1993frontal, 
budd2013monge, visram2014framework, visram2014asymptotic}. 
From a mathematical perspective, the connections with optimal transportation 
have been exposed, leading to new numerical methods and analytical insight 
(see \citet{cullen2007modelling} for a review), whilst 
\citet{cotter2013variational} considered the geometric structure and 
conservation laws of this slice model. 
\citet{visram2014framework} presented a framework for evaluating model 
error in terms of asymptotic convergence in the Eady model. The framework 
is based on the semi-geostrophic (SG) theory in which hydrostatic balance 
and geostrophic balance of the out-of-slice component of the wind are 
imposed to the equations. 
The SG equation provides a suitable limit for asymptotic convergence in 
the Eady model as the Rossby number decreases to zero 
\citep[][]{cullen2008comparison}. 
We use this framework to validate the numerical implementation and assess 
the long term performance of the model developed using the compatible 
finite element method.

The main goal of this paper is to demonstrate the compatible finite element 
approach for NWP in the context of this frontogenesis 
test case.
A challenge in using the compatible finite element method in NWP 
models is the implementation of (the finite element 
version of) the Charney-Philips grid staggering in the vertical direction that 
is used in many current operational forecasting models, such as the Met 
Office Unified Model \citep[][]{wood2014inherently}. 
This requires the temperature space to be a tensor product of discontinuous 
functions in the horizontal direction and continuous functions in the vertical 
direction \citep[][]{cotter2016embedded}. 
Therefore we propose a new advection scheme for a partially continuous finite 
element space and use it to discretise the temperature equation. The other key 
features of the model are: 
i) an upwind DG method is applied for the momentum 
equations; ii) the semi-implicit equations are solved together with an explicit 
strong-stability-preserving Runge-Kutta (SSPRK) scheme to all of the 
advection terms; iii) a balanced initialisation is 
introduced to enforce hydrostatic and geostrophic balances in the initial fields. 

The rest of the paper is structured as follows. 
Section \ref{model_description} provides the model description, including 
formation of the Eady 
problem, discretisation of the governing equations in time and space, and 
settings of the frontogenesis experiment. 
In section \ref{results}, we present the results of the frontogenesis experiments using 
the developed model. 
Here we evaluate model error in terms of SG limit analysis as well as 
energy dynamics. Finally, in section \ref{conclusion} we provide a summary and outlook.

\section{The incompressible Euler-Boussinesq Eady slice model}\label{model_description}
\subsection{Governing equations}

In this section, we describe the model equations for the vertical
slice Eady problem. The equations are as described in
\citet{visram2014framework}, but we repeat them here to
establish notation. To derive a set of equations for the vertical
slice Eady problem, we start from the three-dimensional incompressible
Euler-Boussinesq equations with rigid-lid conditions on the upper and 
lower boundaries,
\begin{eqnarray}
\frac{\partial \MM{u}}{\partial t} + (\MM{u} \cdot \nabla)\MM{u} + f \MM{\hat{z}} 
\times \MM{u} &=& -\frac{1}{\rho_{0}} \nabla p + \frac{g}{\theta_{0}}\theta \MM{\hat{z}},\\
\frac{\partial \theta}{\partial t} + (\MM{u} \cdot \nabla)\theta &=& 0,\\
\nabla \cdot \MM{u} &=& 0,
\end{eqnarray}
where $\MM{u} = (u, v, w)$ is the velocity vector, $\nabla =
(\partial_x, \partial_y, \partial_z)$ is the gradient operator, and
$\MM{\hat{z}}$ is a unit vector in the $z$-direction; $p$ is the
pressure, $\theta$ is the potential temperature, and $g$ is the
acceleration due to gravity; $\rho_{0}$ and $\theta_{0}$ are reference
density and potential temperature values at the surface,
respectively. The rotation frequency $f$ is constant.

In the Eady slice model, we consider perturbations to the constant
background temperature and pressure profiles, $\bar{\theta}(y,z)$ and
$\bar{p}(y,z)$, respectively. All perturbation variables are then
assumed to be independent of $y$, denoted with primed variables as follows,
\begin{eqnarray}
\theta &=& \bar{\theta} (y, z) + \theta'(x, z, t), \\
p &=& \bar{p}(y,z) + p'(x, z, t).
\end{eqnarray}
Following \citet{snyder1993frontal}, we select a background profile 
assuming the geostrophic balance:
\begin{eqnarray}
\bar{\theta}(y, z) &=& \frac{\theta_{0}}{g}(-f \Lambda y + N^{2}z), \label{background_theta}\\
-\frac{1}{\rho_{0}}\frac{\partial \bar{p}}{\partial y} &=& 
f\Lambda \left(z-\frac{H}{2}\right),
\end{eqnarray}
where $\Lambda$ is the constant vertical shear, $N$ is the Brunt-V\"ais\"al\"a frequency, 
and $H$ is the height of the domain. 
The variation in $y$ of the background pressure is therefore written as
\begin{eqnarray}
\frac{\partial \bar{p}}{\partial y} &=& \frac{\rho_{0}g}{\theta_{0}}
\frac{\partial \bar{\theta}}{\partial y}\left(z-\frac{H}{2}\right), \label{dpdy}
\end{eqnarray}
where
\begin{eqnarray}
\frac{\partial \bar{\theta}}{\partial y} &=& -\frac{\theta_{0}f\Lambda}{g} = \mbox{const}.
\end{eqnarray}
Similarly, the variation in $z$ of the background pressure is given by the hydrostatic balance as
\begin{eqnarray}
\frac{\partial \bar{p}}{\partial z} &=& \frac{\rho_{0}g}{\theta_{0}}\bar{\theta}.\label{dpdz} 
\end{eqnarray}

Substituting in the background profiles \eqref{dpdy} and \eqref{dpdz} and 
$\partial/\partial y = 0$ for all perturbation variables, we obtain the nonhydrostatic, 
incompressible Euler-Boussinesq Eady equations in the vertical slice 
with rigid-lid conditions on the upper and lower boundaries,
\begin{eqnarray}
\frac{\partial u}{\partial t} + u\frac{\partial u}{\partial x} + w\frac{\partial u}{\partial z} - fv 
&=& -\frac{1}{\rho_{0}}\frac{\del p'}{\del x},\\
\frac{\partial v}{\partial t} + u\frac{\partial v}{\partial x} + w\frac{\partial v}{\partial z} + fu 
&=& -\frac{g}{\theta_{0}}\frac{\del \bar{\theta}} {\del y}\left(z-\frac{H}{2}\right),\\
\frac{\partial w}{\partial t} + u\frac{\partial w}{\partial x} + w\frac{\partial w}{\partial z} 
&=& - \frac{1}{\rho_{0}}\frac{\del p'}{\del z} + \frac{g}{\theta_0}\theta' ,\\
\frac{\partial \theta'}{\partial t} + u\frac{\partial \theta'}{\partial x} + w\frac{\partial \theta'}{\partial z} 
+ v\frac{\del \bar{\theta}}{\del y} + w\frac{\del \bar{\theta}}{\del z} &=& 0,\\
\frac{\del u}{\del x} + \frac{\del w}{\del z} &=& 0,
\end{eqnarray}
where all variables $u, v, w, \theta'$ and $p'$ are functions of $(x, z, t)$.
Now, we redefine the velocity vector and the gradient operator as those in the vertical slice, 
\begin{eqnarray}
\MM{u} = (u, w), \ \ \ \ \nabla = (\partial_x, \partial_z), 
\end{eqnarray}
and introduce the in-slice buoyancy and the background buoyancy,
\begin{eqnarray}
b' = \frac{g}{\theta_{0}}\theta', \ \ \ \ \bar{b}=\frac{g}{\theta_{0}}\bar{\theta}, \label{def_b}
\end{eqnarray}
respectively. 
Finally, dropping the primes gives the vector form of the Eady slice model equations as
\begin{eqnarray}
\frac{\partial \MM{u}}{\partial t} + (\MM{u} \cdot \nabla)\MM{u} - fv\hat{\MM{x}} 
&=& -\frac{1}{\rho_{0}}\nabla p + b\hat{\MM{z}}, \label{ueq}\\
\frac{\partial v}{\partial t} +  \MM{u} \cdot \nabla v + f\MM{u} \cdot \hat{\MM{x}} 
&=& -\frac{\del \bar{b}} {\del y}\left(z-\frac{H}{2}\right),\label{veq}\\
\frac{\partial b}{\partial t} + \MM{u} \cdot \nabla b + \frac{\partial \bar{b}}{\partial y} v 
+ N^{2} \MM{u} \cdot \hat{\MM{z}} &=& 0,\label{beq}\\
\nabla \cdot \MM{u} &=& 0 \label{peq},
\end{eqnarray}
where $\MM{\hat{x}}$ is a unit vector in the $x$-direction, 
and in \eqref{beq} we used the relationship
\begin{eqnarray}
N^{2} = \frac{g}{\theta_{0}}\frac{\partial \bar{\theta}}{\partial z} = \frac{\partial \bar{b}}{\partial z},
\end{eqnarray}
obtained from \eqref{background_theta} and \eqref{def_b}.

The solutions to the slice model equations \eqref{ueq} to \eqref{peq} are equivalent to 
a $y$-independent solution of the full three dimensional equations. 
The slice model conserves the total energy
\begin{eqnarray}
E = K_{u} + K_{v} + P, \label{total_energy}
\end{eqnarray}
where 
\begin{eqnarray}
K_{u} &=& \rho_{0} \int_{\Omega} \frac{1}{2} |\MM{u}|^{2} \ \diff x, \label{kinetic_u}\\
K_{v} &=& \rho_{0} \int_{\Omega} \frac{1}{2} v^{2} \ \diff x,\label{kinetic_v}\\
P &=& -\rho_{0} \int_{\Omega} b \left(z-\frac{H}{2}\right) \ \diff x, \label{potential}
\end{eqnarray}
are the kinetic energy from the in-slice velocity components, kinetic energy from the 
out-of-slice velocity component,  and potential energy, respectively.

\subsection{Finite element discretisation}
\subsubsection{Compatible finite element spaces}
In this study, a compatible finite element method is used to discretise the 
governing equations. 
First we take our computational domain, denoted by $\Omega$, to be a rectangle 
in the vertical plane with a periodic boundary condition in the $x$-direction,
and rigid-lid conditions on the upper and lower boundaries. 
We refer to the combination of the upper and lower boundaries as $\partial\Omega$.

\begin{figure}[htbp]
  \centering
  \includegraphics[height=65mm]{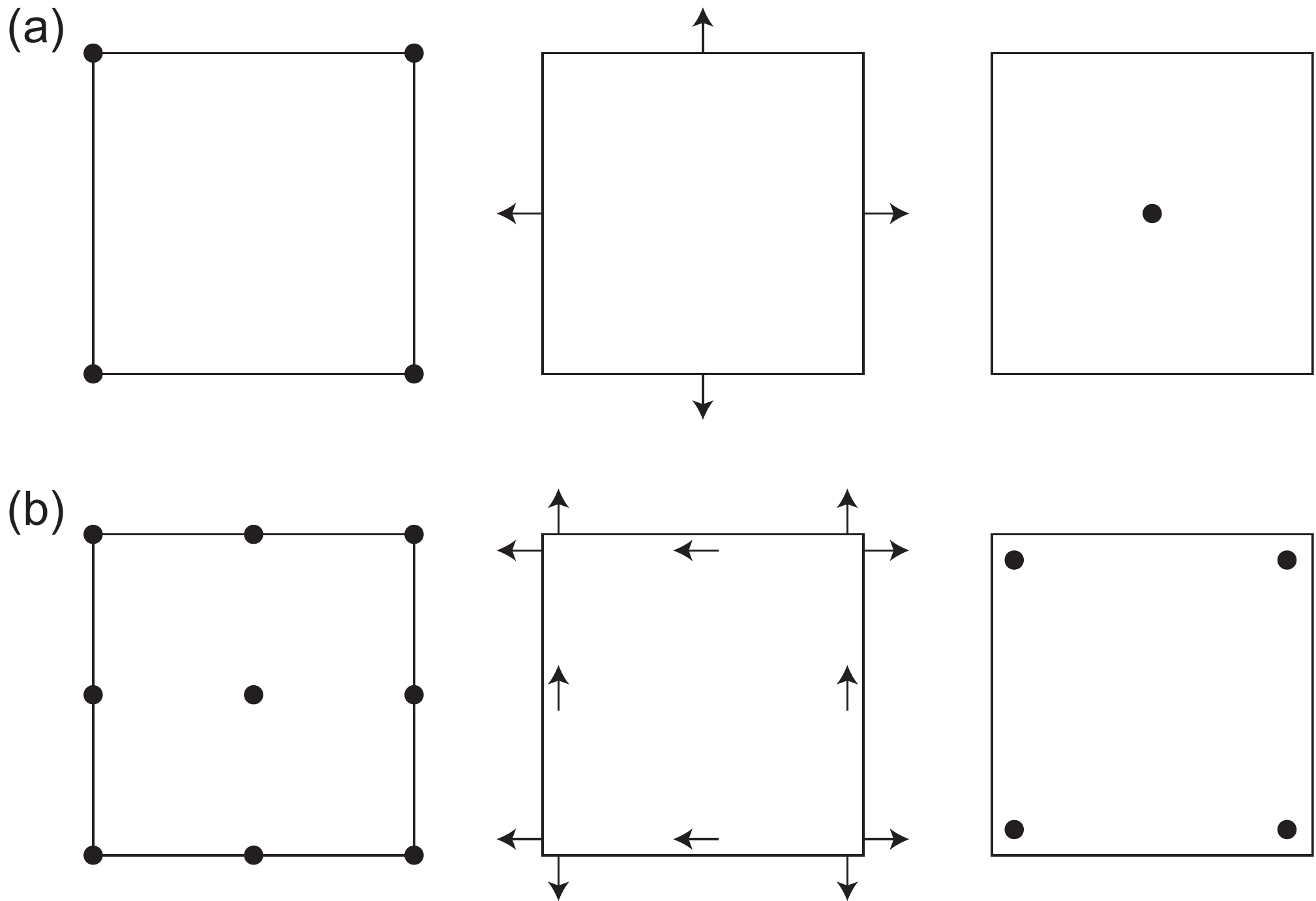}
  \caption{Diagrams showing the nodes for the finite element spaces 
  ($ \mathbb{V}_0$, $ \mathbb{V}_1$, $ \mathbb{V}_2$) 
  = (CG$_k$, RT$_{k-1}$, DG$_{k-1}$) 
  on quadrilaterals in the vertical plane.
  Circles denote scalar nodes, whilst arrows denote normal and tangential 
  components of a vector. Normal components are continuous across element 
  boundaries. Since tangential components are not required to be continuous, 
  these values are not shared by neighbouring elements.
  (a) From left to right: $ \mathbb{V}_0$, $ \mathbb{V}_1$ and $ \mathbb{V}_2$ with $k$ = 1. 
  (b) From left to right: $ \mathbb{V}_0$, $ \mathbb{V}_1$ and $ \mathbb{V}_2$ with $k$ = 2.}
  \label{fig:space-nodes}
\end{figure}

Next we choose finite element spaces with the following properties:
\begin{equation}
\begin{CD}
\underbrace{\mathbb{V}_0(\Omega)}_{\mbox{\small Continuous}} 
@> \nabla^\perp 
>> \underbrace{ \mathbb{V}_1(\Omega)}_{\mbox{\small Continuous\ normal\ components}} 
@>\nabla\cdot >> \underbrace{ \mathbb{V}_2(\Omega)}_{\mbox{\small Discontinuous}},
\end{CD}
\label{spaces}
\end{equation}
where $\nabla^\perp = (-\partial_{z}, \partial_{x})$, $ \mathbb{V}_0$ contains 
scalar-valued continuous functions, $ \mathbb{V}_1$ contains vector-valued 
functions with continuous normal components across element boundaries, 
and $ \mathbb{V}_2$ contains scalar-valued functions that are discontinuous 
across element boundaries. 
The use of these spaces may be considered as an extension of the Arakawa-C 
horizontal grid staggering in finite difference methods.
On quadrilateral elements, \citet{cotter2012mixed} advocated the choice
($ \mathbb{V}_0$, $ \mathbb{V}_1$, $ \mathbb{V}_2$) 
= (CG$_k$, RT$_{k-1}$, DG$_{k-1}$) for $k > 0$, 
where CG$_k$ denotes the continuous finite element space of polynomial degree $k$, 
RT$_k$ denotes the quadrilateral Raviart-Thomas space of polynomial degree $k$, 
and DG$_k$ denotes the discontinuous finite element space of polynomial degree $k$.  
This set of spaces ensures the ratio between global velocity DoFs and 
global pressure DoFs to be exactly 2:1, which helps to avoid spurious modes.
Figure~\ref{fig:space-nodes} provides diagrams in the vertical plane showing the nodes 
for the three spaces for the cases $k = 1$ and $k = 2$.  

We then restrict the model variables to suitable function spaces. First the velocity and 
pressure are defined as
\begin{eqnarray}
\quad \MM{u}\in  \mathring{\mathbb{V}}_1, \quad p \in  \mathbb{V}_2,
\end{eqnarray}
where $\mathring{\mathbb{V}}_1$ is the subspace defined by
\begin{equation}
\mathring{\mathbb{V}}_1 = \left\{
\MM{u}\in \mathbb{V}_1:\MM{u}\cdot\MM{n}=0 \mbox{ on }\partial\Omega\right\}.
\end{equation}
To be consistent with the three-dimensional Arakawa-C grid staggering, we can 
choose $v$ from the same space as $p$ as 
\begin{eqnarray}
\quad v \in  \mathbb{V}_2. 
\end{eqnarray}

\begin{figure}[htbp]
  \centering
  \includegraphics[height=65mm]{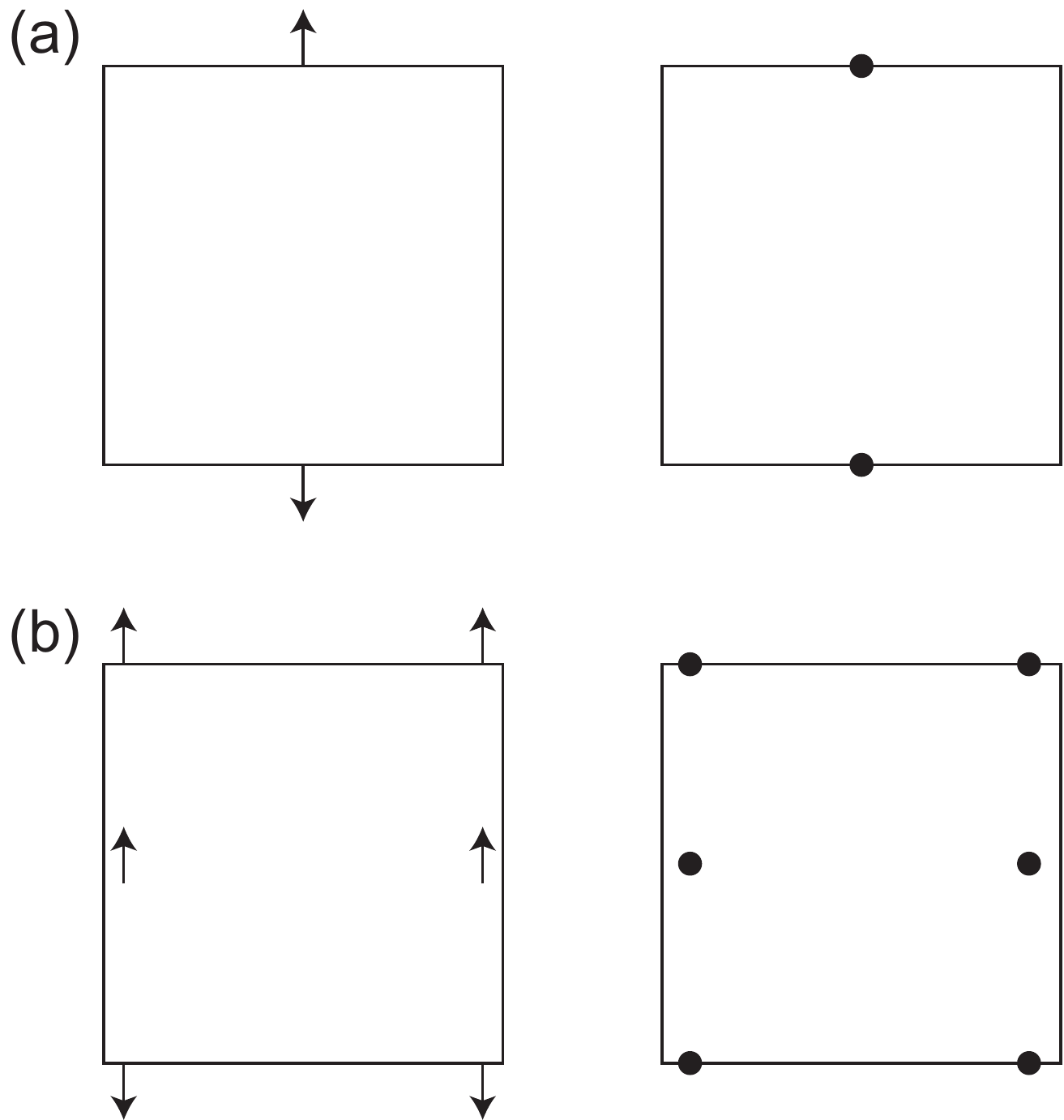}
  \caption{Diagrams showing the nodes for (a) the vertical part of $ \mathbb{V}_1$ (left) 
  and $ \mathbb{V}_b$ (right) with $k$ = 1, and (b) those with $k$ = 2. 
  Circles denote scalar nodes, whilst arrows denote normal and tangential 
  components of a vector.}
  \label{fig:vertical-nodes}
\end{figure}

There are two main options for arranging the temperature/buoyancy in finite difference 
models: the Lorenz grid (temperature collocated with pressure), and the Charney-Phillips 
grid (temperature collocated with vertical velocity). 
To mimic the Lorenz grid, we can simply choose $b$ from the pressure space $\mathbb{V}_2$. 
In this study, we use the Charney-Phillips grid since it avoids spurious hydrostatic 
pressure modes. 
For this purpose, we introduce a scalar space $\mathbb{V}_b$ which is 
obtained by the tensor product of the DG$_{k-1}$ space in the horizontal direction 
and the CG$_{k}$ space in the vertical direction. 
As shown in Figure~\ref{fig:vertical-nodes},
$\mathbb{V}_b$ has the same node locations as the vertical part of $ \mathbb{V}_1$.
Then we choose the buoyancy as
\begin{eqnarray}
\quad b \in  \mathbb{V}_b. 
\end{eqnarray}
This constructs the extension of the Charney-Phillips staggering to
compatible finite element spaces. 
\citet{natale2016compatible} showed that this choice of finite element spaces leads to
a one-to-one mapping between pressure and buoyancy in the hydrostatic balance equation
and therefore, as in the finite difference models using the Charney-Phillips grid,
avoids spurious hydrostatic pressure modes.
Since the space $ \mathbb{V}_b$ is discontinuous in the horizontal direction and 
continuous in the vertical direction, we need a transport scheme for a partially 
continuous finite element space, which we detail in the next subsection.

\subsubsection{Spatial discretisation}\label{space}
We now use the compatible finite element spaces introduced above to discretise the 
model equations \eqref{ueq} to \eqref{peq}. 
Here we start with the discretisation of the in-slice velocity equation \eqref{ueq}. 
First we rewrite the advection term as
\begin{equation}
(\MM{u}\cdot\nabla)\MM{u} 
= (\nabla^{\perp}\cdot\MM{u})\MM{u}^{\perp}+\frac{1}{2}\nabla|\MM{u}|^{2}. \label{invariant}
\end{equation}
Then, taking \eqref{ueq}, dotting with a test function 
$\MM{w} \in \mathring{\mathbb{V}}_1$, 
and integrating over the domain gives
\begin{eqnarray}
\int_\Omega\MM{w}\cdot \frac{\partial \MM{u}}{\partial t} \diff x 
+ \int_\Omega\MM{w}\cdot(\nabla^{\perp}\cdot\MM{u})\MM{u}^{\perp}\diff x 
= \int_\Omega\MM{w}\cdot fv\MM{\hat{x}} \diff x  - \int_\Omega \MM{w} \cdot 
\nabla \left(\frac{p}{\rho_{0}} + \frac{1}{2}|\MM{u}|^{2}\right) \diff x + \int_\Omega 
\MM{w}\cdot b\hat{\MM{z}} \diff x, \ \ \forall\MM{w}\in  \mathring{\mathbb{V}}_1. 
\label{udis} 
\end{eqnarray}
Recalling that $\MM{u}$ is in $ \mathring{\mathbb{V}}_1$, $\nabla^{\perp}\cdot\MM{u}$ 
in the second term is not generally defined, since the tangential component of $\MM{u}$ 
is not continuous across element boundaries in general. 
We resolve this by integrating the term by parts.
For the contribution to the integral from each element $e$ we obtain
\begin{equation}
\int_e \MM{w}\cdot(\nabla^\perp\cdot\MM{u})\MM{u}^\perp \diff x
= -\int_e \nabla^\perp(\MM{w}\cdot\MM{u}^\perp)\cdot\MM{u} \diff x
    +\int_{\partial e}\MM{n}^\perp\cdot\tilde{\MM{u}}\MM{w}\cdot\MM{u}^\perp \diff S,
\end{equation}
where $\tilde{\MM{u}}$ is the upwind value of $\MM{u}$ on the element boundary
$\partial e$.
Summing over all elements, the advection term becomes
\begin{equation}
\int_\Omega \MM{w}\cdot(\nabla^\perp\cdot\MM{u})\MM{u}^\perp \diff x = 
-\int_\Omega \nabla^\perp(\MM{w}\cdot\MM{u}^\perp)\cdot\MM{u} \diff x 
+\int_\Gamma \jump{\MM{w}\cdot\MM{u}^\perp}^{\perp}\cdot\tilde{\MM{u}} \diff S, 
\label{uadv_int}
\end{equation}
where $\Gamma$ is the set of interior facets in the finite element mesh with the 
two sides of each facet arbitrarily labelled by + and −, the jump operator is defined by
\begin{eqnarray}
\jump{q} &=& q^{+}\MM{n}^{+} + q^{-}\MM{n}^{-}, \\
\jump{\MM{v}} &=& \MM{v}^{+}\cdot\MM{n}^{+} + \MM{v}^{-}\cdot\MM{n}^{-},
\end{eqnarray}
for any scalar $q$ and vector $\MM{v}$, and $\tilde{\MM{u}}$ is evaluated on the 
upwind side as
\begin{eqnarray}
\tilde{\MM{u}} = \left\{
\begin{array}{l}
\MM{u}^{+}  \ \ \ \ \mbox{if} \ \MM{u}\cdot\MM{n}^{+} < 0, \\
\MM{u}^{-} \ \ \ \ \mbox{otherwise.}
\end{array}
\right.
\end{eqnarray}
A variational derivation and numerical analysis of the discretisation of this
term can be found in \citet{natale2016variational}.
Turning attention to the pressure gradient term $\nabla \left(\frac{p}{\rho_{0}} 
+ \frac{1}{2}|\MM{u}|^{2}\right)$, we also integrate by parts.
The discrete form of the in-slice velocity equation becomes
\begin{eqnarray}
\int_{\Omega}\MM{w}\cdot \frac{\partial \MM{u}}{\partial t} \diff x 
&=&  \int_\Omega \nabla^\perp(\MM{w}\cdot \MM{u}^{\perp})\cdot\MM{u} \diff x 
+ \int_{\Omega} \nabla\cdot\MM{w} \left(\frac{p}{\rho_{0}} 
+ \frac{1}{2}|\MM{u}|^{2}\right) \diff x \nonumber \\
&+& \int_{\Omega}\MM{w}\cdot fv\MM{\hat{x}} \diff x  
+ \int_{\Omega} \MM{w}\cdot b\hat{\MM{z}} \diff x 
- \int_\Gamma \jump{\MM{w}\cdot\MM{u}^\perp}^\perp\cdot\tilde{\MM{u}} \diff S, 
\ \ \forall\MM{w}\in  \mathring{\mathbb{V}}_1. \label{ueq_int} 
\end{eqnarray}

The out-of-slice velocity space $\mathbb{V}_2$ is discontinuous. An upwind DG treatment 
of \eqref{veq} leads to
\begin{eqnarray}
\int_{\Omega}\phi\frac{\partial v}{\partial t} \diff x 
-  \int_\Omega \nabla \cdot (\phi \MM{u}) v \diff x 
+ \int_\Gamma \jump{\phi \MM{u}} \tilde{v} \diff S 
+ \int_{\Omega} \phi f \MM{u} \cdot \hat{\MM{x}} \diff x 
+ \int_{\Omega} \phi \frac{\del \bar{b}} {\del y} \left(z-\frac{H}{2}\right) \diff x = 0, 
\ \ \forall\phi\in  \mathbb{V}_2,  \label{veq_int} 
\end{eqnarray}
where $\tilde{v}$ denotes the upwind value of $v$.

We now describe how we discretise the buoyancy equation \eqref{beq}. 
Recall that $b$ is in the finite element space $ \mathbb{V}_b$, which is obtained 
by the tensor product of a discrete finite element space in the horizontal direction 
with a continuous finite element space in the vertical direction. 
Therefore we propose a blend of an upwind DG method in the horizontal direction 
and Streamline Upwind Petrov Galerkin (SUPG) method in the vertical direction. 
First, multiplying the equation \eqref{beq} by a test function $\gamma \in  \mathbb{V}_b$ 
and integrating it over each column $C$ gives
\begin{eqnarray}
\int_{C} \gamma \frac{\partial b}{\partial t} \diff x 
+ \int_{C} \gamma \MM{u} \cdot \nabla b \diff x 
+ \int_{C} \gamma \frac{\del \bar{b}}{\del y} v \diff x 
+ \int_{C} \gamma N^{2}w \diff x = 0, \ \ \forall\gamma\in  \mathbb{V}_b.
\end{eqnarray}
To obtain the DG formulation, we apply integration by parts in each column,
\begin{eqnarray}
\int_{C} \gamma \frac{\partial b}{\partial t} \diff x 
- \int_{C} \nabla \cdot (\gamma \MM{u}) b \diff x 
+  \int_{C} \gamma \frac{\del \bar{b}}{\del y} v \diff x 
+ \int_{C} \gamma N^{2} w \diff x 
+ \int_{\partial C} \gamma \MM{u} \cdot \MM{n} \tilde{b} \diff S = 0, 
\ \ \forall\gamma\in  \mathbb{V}_b,
\end{eqnarray}
where $\tilde{b}$ denotes the upwind value of $b$ on the column boundary $\partial C$. 
Integrating by parts again gives
\begin{eqnarray}
\int_{C} \gamma \frac{\partial b}{\partial t} \diff x 
+ \int_{C} \gamma \MM{u} \cdot \nabla b \diff x 
+  \int_{C} \gamma \frac{\del \bar{b}}{\del y} v \diff x + \int_{C} \gamma N^{2} w \diff x 
+ \int_{\partial C} \gamma \MM{u} \cdot \MM{n} \tilde{b} 
- \gamma \MM{u} \cdot \MM{n} b \diff S = 0, \ \ \forall\gamma\in  \mathbb{V}_b,
\end{eqnarray}
where the new boundary term contains the value of $b$ on the interior of the column. 
Summing over the whole domain, we obtain
\begin{eqnarray}
\int_{\Omega} \gamma \frac{\partial b}{\partial t} \diff x 
+ \int_{\Omega} \gamma \MM{u} \cdot \nabla b \diff x 
+  \int_{\Omega} \gamma \frac{\del \bar{b}}{\del y} v \diff x 
+ \int_{\Omega} \gamma N^{2} w \diff x + \int_{\Gamma_v} \jump{\gamma \MM{u}} \tilde{b} 
- \jump{\gamma \MM{u} b} \diff S = 0, \ \ \forall\gamma\in  \mathbb{V}_b,
\end{eqnarray}
where $\Gamma_v$ denotes the set of vertical facets of each column. 
We then apply the SUPG method in the vertical direction by replacing the test function 
$\gamma$ with $\gamma + \tau \gamma_z$, where $\gamma_z$ denotes the vertical
derivative of $\gamma$:
\begin{eqnarray}
\int_{\Omega} (\gamma + \tau \gamma_z) \frac{\partial b}{\partial t} \diff x 
+ \int_{\Omega} (\gamma + \tau \gamma_z) \MM{u} \cdot \nabla b \diff x 
+ \int_{\Omega} \gamma \frac{\del \bar{b}}{\del y} v \diff x 
+  \int_{\Omega} (\gamma + \tau  \gamma_z) N^{2} w \diff x \nonumber \\
+ \int_{\Gamma_v} \jump{(\gamma + \tau \gamma_z) \MM{u}} \tilde{b} 
- \jump{(\gamma + \tau \gamma_z) \MM{u} b} \diff S = 0, 
\ \ \forall\gamma\in  \mathbb{V}_b. \label{beq_int}
\end{eqnarray}
Here $\tau$ is an upwinding coefficient
\begin{eqnarray}
\tau = c \Delta t \MM{u} \cdot \hat{\MM{z}},
\end{eqnarray}
where $\Delta t$ is the time step and the constant $c$ is set at 
$1/15^{\frac{1}{2}}$ following \citet{raymond1976selective}.

Finally we discretise the continuity equation \eqref{peq}. 
Multiplying by a test function $ \sigma \in  \mathbb{V}_2$ and integrating it 
over the domain $\Omega$ gives 
\begin{eqnarray}
\int_\Omega \sigma \nabla \cdot \MM{u} \diff x = 0, \forall \sigma\in  \mathbb{V}_2. 
\label{cont_int}
\end{eqnarray}
Since $\nabla \cdot \MM{u}$ can be defined globally in $ \mathbb{V}_2$, 
the projection of \eqref{cont_int} is trivial, i.e., the incompressible condition is 
satisfied exactly under this discretisation.

\subsubsection{Time discretisation}\label{time_discretisation}
Now we discretise the equations \eqref{ueq_int}, \eqref{veq_int} and 
\eqref{beq_int} in time using a semi-implicit time-discretisation scheme. 
This is most easily described as a fixed number of iterations for a Picard 
iteration scheme applied to a fully implicit time integration scheme.

The implicit time integration scheme is obtained by applying a
(possibly off-centred) implicit time discretisation average to all of
the forcing terms, as well as the advecting velocity in all of the
equations. We then apply an explicit SSPRK scheme to all of the
advection terms. To write down the scheme, we first define operators 
$L_{\MM {u}}$, $L_{v}$ and $L_{b}$ in the following implicit time-stepping 
formulation,
\begin{eqnarray}
\int_{\Omega}\MM{w}\cdot L_{\MM {u}}\MM{u} \diff x 
&=& \Delta t \int_\Omega \nabla^\perp(\MM{w}\cdot\MM{u}^{*\perp})\cdot\MM{u} \diff x 
+ \Delta t \int_{\Omega} \nabla\cdot\MM{w} \left(\frac{p^*}{\rho_{0}} 
+ \frac{1}{2}|\MM{u}^*|^{2}\right) \diff x 
+ \Delta t  \int_{\Omega}\MM{w}\cdot fv^*\MM{\hat{x}} \diff x \nonumber \\
 &+& \Delta t \int_{\Omega} \MM{w}\cdot b^*\hat{\MM{z}} \diff x 
 - \Delta t \int_\Gamma \jump{\MM{w}\cdot\MM{u}^{*\perp}}^\perp\cdot\tilde{\MM{u}}^* \diff S, 
 \ \ \forall\MM{w}\in  \mathring{\mathbb{V}}_1, \label{ueq_L} \\
\int_{\Omega}\phi L_{v} v \diff x 
&=& \Delta t \int_\Omega \nabla \cdot (\phi \MM{u}^*) v \diff x 
- \Delta t \int_{\Omega} \phi f \MM{u}^* \cdot \hat{\MM{x}} \diff x 
- \Delta t \int_{\Omega} \phi \frac{\del \bar{b}} {\del y} \left(z-\frac{H}{2}\right) \diff x 
\nonumber \\
 &-& \Delta t \int_\Gamma \jump{\phi \MM{u}^*} \tilde{v}^* \diff S , 
 \ \ \forall\phi\in  \mathbb{V}_2,  \label{veq_L} \\
\int_{\Omega} (\gamma + \tau \gamma_z) L_{b} b \diff x 
&=& - \Delta t \int_{\Omega} (\gamma + \tau \gamma_z) \MM{u}^* \cdot \nabla b \diff x 
- \Delta t \int_{\Omega} (\gamma + \tau \gamma_z) \frac{\del \bar{b}}{\del y} v^* \diff x
 - \Delta t \int_{\Omega} (\gamma + \tau \gamma_z) N^{2} w^* \diff x \nonumber \\
&-& \Delta t \int_{\Gamma_v} \jump{(\gamma + \tau \gamma_z) \MM{u}^*} \tilde{b}^* 
- \jump{(\gamma + \tau \gamma_z) \MM{u}^* b^*} \diff S, 
\ \ \forall\gamma\in  \mathbb{V}_b, \label{beq_L} 
\end{eqnarray}
where the star denotes $y^{*} = (1-\alpha)y^{n} + \alpha y^{n+1}$ with a
time-centring parameter $\alpha$. 
A 3rd order 3 step SSPRK time-stepping method \citep{shu1988efficient} 
is then applied as
\begin{eqnarray}
\varphi^{1}_{y} &=& y^{n}+L_{y}y^{n}, \\
\varphi^{2}_{y} &=& \frac{3}{4}y^{n}+\frac{1}{4}(\varphi^{1}_{y}+L_{y} \varphi^{1}_{y}), \\
Ay^{n} &=& \frac{1}{3}y^{n}+\frac{2}{3}(\varphi^{2}_{y}+L_{y} \varphi^{2}_{y}), \label{eq:Advection_step}
\end{eqnarray}
for each variable $y = \MM{u}, v$ and $b$, where $A$ is the advection operator.

Finally, we solve for $\MM{u}^{n+1}$, $v^{n+1}$, $b^{n+1}$ and $p^{n+1}$ 
iteratively using a Picard iteration method,
\begin{eqnarray}
\int_{\Omega}\MM{w}\cdot \Delta \MM{u} \diff x - \alpha \Delta t \int_{\Omega} 
\nabla\cdot\MM{w} \left(\frac{\Delta p}{\rho_{0}} \right) \diff x 
- \alpha \Delta t \int_{\Omega}\MM{w}\cdot f \Delta v \MM{\hat{x}} \diff x  
\nonumber \\
- \alpha \Delta t \int_{\Omega} \MM{w}\cdot \Delta b\hat{\MM{z}} \diff x 
&=& -R_u[\MM{w}], \ \ \forall\MM{w}\in  \mathring{\mathbb{V}}_1,   \\
\int_{\Omega}\phi \Delta v \diff x + \alpha \Delta t \int_{\Omega} \phi f \Delta u \diff x 
&=& -R_v[\phi], \ \ \forall\phi\in  \mathbb{V}_2, \\
\int_{\Omega} \gamma \Delta b \diff x 
+ \alpha \Delta t \int_{\Omega} \gamma N^{2} \Delta w \diff x 
&=& -R_b[\gamma], \ \ \forall\gamma\in  \mathbb{V}_b, \\
\int_{\Omega}  \sigma \nabla \cdot \Delta \MM{u}  \diff x  
&=& -R_p[ \sigma], \ \ \forall \sigma\in  \mathbb{V}_2.
\end{eqnarray}
Here $R_{\MM{u}}[\MM{w}]$, $R_v[\phi]$, $R_b[\gamma]$ and $R_p[ \sigma]$ 
are the residuals for the implicit system,
\begin{eqnarray}
R_u[\MM{w}] &=& \int_{\Omega} (\MM{u}^{n+1} - A\MM{u}^{n})\cdot\MM{w} \diff x, 
\ \ \forall\MM{w}\in  \mathring{\mathbb{V}}_1, 
\label{eq:R_u_w}\\
R_v[\phi] &=& \int_{\Omega} (v^{n+1} - Av^{n})\phi \diff x, \ \ \forall\phi\in  \mathbb{V}_2, 
\label{eq:R_v_phi}\\
R_b[\gamma] &=& \int_{\Omega} (b^{n+1} - Ab^{n})\gamma \diff x,  
\ \ \forall\gamma\in  \mathbb{V}_b, \label{eq:R_b_gamma} \\
  R_p[\sigma] &=& \int_{\Omega} \nabla \cdot \MM{u}^{n+1} \, \sigma \diff x, 
  \ \ \forall \sigma\in  \mathbb{V}_2, \label{eq:R_p_sigma}
\end{eqnarray}
where $A$ is as defined in \eqref{eq:Advection_step}.
After obtaining $\Delta \MM{u}$, $\Delta v$, $\Delta b$ and $\Delta p$, 
we replace $\MM{u}^{n+1}$, $v^{n+1}$, $b^{n+1}$ and $p^{n+1}$
with $\MM{u}^{n+1} + \Delta \MM{u}$ , $v^{n+1} + \Delta v$ , 
$b^{n+1} + \Delta b$ and $p^{n+1} + \Delta p$, respectively. Then we repeat the
iterative procedure for a fixed number of times, which is set to 4 in this study. 
Figure \ref{fig:pseudocode} provides pseudocode for the timestepping procedure.
As the 3rd order SSPRK schemes are stable for both DG and SUPG methods, the
system is well conditioned for stable Courant numbers, and can be
solved with a few iterations of preconditioned GMRES applied to the
full coupled system of 4 variables.

\begin{figure}[t]
    \begin{algorithmic}
      \State $y^{n}, \; y^{n+1} \gets y^0$
      \Comment initialise
      \For{$k < k_{\mathrm{max}}$}
      \Comment time step loop
      \For{$i < i_{\mathrm{max}}$}
      \Comment Picard iteration loop
      \State \textbf{Compute:} $Ay^n$      
      \Comment advection step using (51)--(53)
      \State $R_y \gets y^{n+1}-Ay^n$      
      \Comment update residual 
      \State \textbf{Solve:} $J[\Delta y] = -R_y$
      \Comment linear solve using (54) -- (57)
      \State $y^{n+1} \gets y^{n+1} + \Delta y$ 
      \Comment update variables
      \EndFor
      \State $y^{n} \gets y^{n+1}$ 
      \Comment advance one time step
      \EndFor
    \end{algorithmic}
  \caption{Pseudocode for the timestepping procedure. The variable $y$ represents 
    the model variables $\MM{u}, v, b$ and $p$, and $J[\Delta y]$ denotes 
    the Jacobian from the linear system.
    The constant $i_\mathrm{max}$ denotes a fixed number for a Picard iteration, 
    and $k_\mathrm{max}$ denotes the total number of time steps.}
  \label{fig:pseudocode}
\end{figure}

We use a block diagonal ``Riesz-map'' preconditioner
\citep{mardal2011preconditioning} for the GMRES iterations.  This
operator has an $H(\text{div})$ inner product in the velocity block
and mass matrices in the other diagonal blocks; see
\citet{natale2016compatible} for details of $H(\text{div})$ finite
element spaces.  Inverting the mass matrices is straightforward,
requiring only a few iterations of a stationary iteration such as
Jacobi; the $H(\text{div})$ block is more challenging due to the
non-trivial kernel.  We use an LU factorisation, provided
by MUMPS -- the MUltifrontal Massively Parallel sparse direct Solver
\citep{MUMPS01,MUMPS02} -- to invert it, since we do not currently
have access to a suitable preconditioner for this operator.

\subsection{Experimental settings}
\subsubsection{Constants}\label{constants}

In the frontogenesis experiments, the model constants are set to the values below, 
following \citet{nakamura1989nonlinear}, \citet{cullen2008comparison},
\citet{visram2014framework}, and \citet{visram2014asymptotic}:
\begin{eqnarray}
L &=& 1000\ \mathrm{km},\ H = 10\ \mathrm{km},\ f = 10^{-4}\ \mathrm{s}^{-1},
\ \nonumber\\
g &=& 10\ \mathrm{m\ s}^{-2},\ \rho_{0} = 1\ \mathrm{kg\ m}^{-3},\ \theta_{0} 
= 300\ \mathrm{K}, \nonumber\\
\Lambda &=& 10^{-3}\ \mathrm{s}^{-1},\ N^{2} = 2.5 \times 10^{-5}\ \mathrm{s}^{-1},
 \nonumber
\end{eqnarray}
where $L$ and $H$ determine the model domain $\Omega = [-L, L] \times [0,H]$. 
The $y$ component of the background buoyancy in \eqref{veq} and \eqref{beq} 
is therefore calculated as
\begin{eqnarray}
\frac{\del \bar{b}}{\del y} = -f\Lambda = -10^{-7}\ \mathrm{s}^{-2}. 
\end{eqnarray}
  
The Rossby and Froude numbers are given in the model as
\begin{eqnarray}
\mathrm{Ro} = \frac{u_{0}}{fL} = 0.05, \label{Rossby} \\
\mathrm{Fr} = \frac{u_{0}}{NH} = 0.1,
\end{eqnarray}
where $u_{0} = 5\ \mathrm{m\ s}^{-1}$ is a representative velocity. 
The ratio of Rossby number to the Froude number defines the Burger number,
\begin{eqnarray}
\mathrm{Bu} = \mathrm{Ro}/\mathrm{Fr} = 0.5, 
\end{eqnarray}
which is used when initialising the model.

\subsubsection{Initialisation}\label{initialisation}
The model field is initialised with a small perturbation with the wavelength 
corresponding to the most unstable mode. 
In this study, the following form of the small perturbation is applied to the in-slice buoyancy, 
\begin{eqnarray}
\label{perturbation}
b(x, z) &=& aN \left\{-\left[1-\frac{\mathrm{Bu}}{2}\coth\left(\frac{\mathrm{Bu}}{2}\right)\right] 
\sinh Z \cos \left( \frac{\pi x}{L} \right) -n \mathrm{Bu} \cosh Z \sin 
\left( \frac{\pi x}{L} \right) \right\}, 
\label{init_buoyancy}
\end{eqnarray}
which is the structure of the normal mode taken from \citet{williams1967atmospheric}. 
The constant $a$ corresponds to the amplitude of the perturbation,
and the constant $n$ takes the form of
\begin{eqnarray}
n = \frac{1}{\mathrm{Bu}}\left\{ \left[\frac{\mathrm{Bu}}{2} - 
\tanh \left(\frac{\mathrm{Bu}}{2} \right) \right] 
\left[\coth \left( \frac{\mathrm{Bu}}{2} - \frac{\mathrm{Bu}}{2} \right) \right] \right\}^{\frac{1}{2}}.
\end{eqnarray}
The modified vertical coordinate $Z$ is defined as
 \begin{eqnarray}
Z = \mathrm{Bu}\left[\left(\frac{z}{H} - \frac{1}{2} \right) \right].
\end{eqnarray}

Next we initialise the pressure $p$. Given the $b$ in \eqref{perturbation}, 
we seek a pressure in hydrostatic balance,
\begin{equation}
\label{pvbalance}
\frac{\partial p}{\partial z} = \rho_{0}b.
\end{equation}
Since we have rigid-lid conditions at the upper and lower boundaries, 
we require a symmetry condition for the pressure,
\begin{equation}
\label{sym}
\int_{z=0}^{z=H} p(x,z) \diff z = 0, \quad \forall x.
\end{equation}
This boundary condition is hard to enforce in the solver, 
so we first solve for a hydrostatic pressure $\hat{p}$ with a free surface boundary condition 
on the top, 
\begin{eqnarray}
\frac{\partial \hat{p}}{\partial z} - \rho_{0} b = 0, \quad \hat{p}(z=H) = 0.
\end{eqnarray}
The finite element approximation can be found with a test function 
$\gamma \in  \mathbb{V}_b$ as
\begin{equation}
\int_\Omega \gamma \frac{\partial \hat{p}}{\partial z} \diff x 
- \int_\Omega \gamma \rho_{0} b \diff x 
= - \int_\Omega \frac{\partial\gamma}{\partial z} \hat{p} \diff x 
- \int_\Omega \gamma \rho_{0} b \diff x = 0, 
\ \  \forall\gamma\in  \mathbb{V}_b,
\end{equation}
where we have integrated the pressure gradient term by parts in the second line. 
We then add an arbitrary function of $x$ to the solution $\hat{p}$ to find $p$ 
satisfying the symmetry condition \eqref{sym}. 

Next we initialize the out-of-slice velocity $v$ by seeking a velocity in geostrophic 
balance with the initialised $p$ as
\begin{eqnarray}
\frac{\partial p}{\partial x} = \rho_{0} fv. 
\end{eqnarray}
As $\partial p/\partial x$ is not defined in our finite element framework, 
first we find $\MM{s} = \nabla p$ in $\mathring{\mathbb{V}}_1$ as
\begin{eqnarray}
\int_{\Omega} \MM{w} \cdot \MM{s} \diff x = \int_{\Omega} \MM{w} \cdot \nabla p 
\diff x = -\int_{\Omega} \nabla \cdot \MM{w} \, p \diff x + \int_{\Gamma} \MM{w} 
\cdot \MM{n} p \diff S,  \ \ \forall\MM{w}\in  \mathring{\mathbb{V}}_1,
\end{eqnarray}
where we have integrated the pressure gradient term by parts in the last equality. 
Then we solve for the initial $v$ as
\begin{eqnarray}
\int_{\Omega} \phi \rho_{0} f v \diff x 
= \int_{\Omega} \phi \MM{s}\cdot\hat{\MM{x}} \diff x, 
\ \ \forall\MM{\phi}\in  \mathbb{V}_2. \label{init_v}
\end{eqnarray}

To initialise the in-slice velocity $\MM{u} = (u, w)$, we seek a solution 
to the linear equations for $v$ and $b$,
\begin{eqnarray}
\pp{v_g}{t} &=& -fu - \pp{\bar{b}}{y}
\left(z-\frac{H}{2}\right), \\
\pp{b_g}{t} &=& -\pp{\bar{b}}{y}{v_g} - N^2w, 
\end{eqnarray}
where we make the SG approximation that $v_g$ and $b_g$ are given by the 
geostrophic and hydrostatic balance.  If the pressure is in geostrophic and 
hydrostatic balance then we have
\begin{eqnarray}
\nabla p = \rho_{0}\begin{pmatrix}
fv_g \\
b_g \\
\end{pmatrix},
\end{eqnarray}
and therefore we have
\begin{eqnarray}
\nabla \dot{p} = \rho_{0}
\begin{pmatrix}
f\dot{v}_g \\
\dot{b}_g \\
\end{pmatrix}
= \rho_{0} \begin{pmatrix}
  -f^2 u - f\pp{\bar{b}}{y}
\left(z-\frac{H}{2}\right)\\
-\pp{\bar{b}}{y}{v_g} - N^2w \\
\end{pmatrix},
\end{eqnarray}
where the dot denotes $\dot{y} = \pp{y}{t}$. We rewrite this in a vector form as
\begin{eqnarray}
\rho_{0}
\begin{pmatrix}
f^2 & 0 \\
0 & N^2 \\
\end{pmatrix}
\MM{u} + \nabla \dot{p} = \rho_{0} \pp{\bar{b}}{y}
\begin{pmatrix}
 -f 
\left(z-\frac{H}{2}\right)\\
-{v_g} \\
\end{pmatrix}.
\end{eqnarray}
The finite element approximation is then
\begin{eqnarray}
\int_{\Omega} \MM{w} \cdot \rho_{0}
\begin{pmatrix}
f^2 & 0 \\
0 & N^2 \\
\end{pmatrix}
\MM{u} \diff x - \int_\Omega \nabla \cdot \MM{w} \dot{p} \diff x
= \int_{\Omega} \MM{w}\cdot \rho_{0}\pp{\bar{b}}{y}
\begin{pmatrix}
 -f 
\left(z-\frac{H}{2}\right)\\
-{v_g} \\
\end{pmatrix}\diff x, \quad \forall \MM{w} \in \mathring{\mathbb{V}}_1,
\label{balance_int}
\end{eqnarray}
where we have integrated the pressure gradient term by parts in the first line. 
Here we introduce $\psi \in \mathbb{V}_0$, with $\nabla^\perp \psi=\MM{u}$, 
and choose $\MM{w}=\nabla^\perp \xi$, then we have
\begin{align}
\int_{\Omega} \nabla \xi \cdot 
\begin{pmatrix}
N^2 & 0 \\
0 & f^2 \\
\end{pmatrix}
\nabla \psi  \diff x
= \int_{\Omega} \nabla\xi\cdot \pp{\bar{b}}{y}
\begin{pmatrix}
-{v_g} \\
f\left(z-\frac{H}{2}\right)\\
\end{pmatrix}\diff x, 
\quad \forall \xi \in \mathbb{V}_0.
\end{align} 
With boundary conditions $\psi=0$ on top and bottom, and substituting the 
initialised $v$ from \eqref{init_v} for $v_g$, we obtain the balanced $\psi$. 
We then solve for the initial $\MM{u}$ from $\psi$ to complete the initialisation 
of the model field.

Finally, we introduce a breeding procedure used by \citet{visram2014framework} 
and \citet{visram2014asymptotic} to remove any remaining unbalanced modes 
in the initial condition. In the experiments performed in section \ref{results}, 
the model field is initialised with a small perturbation by choosing $a$ = -7.5 in 
\eqref{init_buoyancy}. The simulation is then advanced for three computational days 
until the maximum amplitude of $v$ reaches 3 m s$^{-1}$, 
at which point the time is reset to zero to match the amplitude of the initial perturbation 
with that of \citet{nakamura1989nonlinear} and \citet{visram2014framework}
as closely as possible. 

\subsubsection{Asymptotic limit analysis}\label{settings_sg}
To validate the numerical implementation and assess the long term performance 
of the model, we introduce the asymptotic limit analysis based on the SG theory. 
The test outlined here follows that described in \citet{cullen2008comparison}, 
and used in \citet{visram2014framework} and \citet{visram2014asymptotic}. 

First we apply the SG approximation to the governing equations \eqref{ueq} to 
\eqref{peq} by imposing the hydrostatic balance and the geostrophic balance 
of the $v$ component of the wind,
\begin{eqnarray}
- fv \hat{\MM{x}} &=& -\frac{1}{\rho_{0}}\nabla p + b\hat{\MM{z}}. \label{ueq_sg}
\end{eqnarray}
The equation \eqref{ueq_sg}, together with the equations \eqref{veq} to \eqref{peq}, 
are the SG equations of the Eady slice model. Now, the solutions of the SG 
equations are invariant to the changes of variables, 
\begin{equation}
x \rightarrow \beta x,\ \ u \rightarrow \beta u,\ \ f \rightarrow \frac{f}{\beta}, \label{rescaling}
\end{equation}
where $\beta$ is a rescaling parameter and all other variables are invariant. 
Recalling the definition of the Rossby number \eqref{Rossby}, the rescaling \eqref{rescaling} 
converts $\mathrm{Ro} \rightarrow \beta \mathrm{Ro}$. 
The SG solution therefore provides an asymptotic limit of the model as the Rossby 
number tends to zero. 

With the rescaling parameter $\beta$, the limit of $\mathrm{Ro} \rightarrow 0$ is 
equivalent to $\beta \rightarrow 0$. 
Thus the convergence of the model to the SG solution can be tested by performing 
a sequence of simulations with decreasing $\beta$. Here we follow 
\citet{cullen2008comparison} in defining the out-of-slice geostrophic imbalance as
\begin{equation}
\eta = v - \frac{1}{\rho_{0}f}\frac{\partial p}{\partial x}, \label{gi}
\end{equation}
which is expected to converge at a rate proportional to $\mathrm{Ro}^{2}$, i.e. $\beta^{2}$. 
To calculate $\eta$ in our finite element framework, first we find the vector of geostrophic 
and hydrostatic imbalance  $\MM{q}$ in $\mathring{\mathbb{V}}_1$ defined as
\begin{eqnarray}
\MM{q} =  \rho_{0} \begin{pmatrix}
fv\\
b\\
\end{pmatrix} - \nabla p.
\end{eqnarray}
where 
\begin{eqnarray}
\eta = \frac{1}{\rho_{0}f} \MM{q}\cdot\MM{x}. \label{gi_v}
\end{eqnarray}
The finite element approximation is obtained as
\begin{eqnarray}
\int_{\Omega} \MM{w} \cdot \MM{q} \diff x =  \int_{\Omega} \MM{w} \cdot \rho_{0} 
\begin{pmatrix}
fv\\
b\\
\end{pmatrix} \diff x - \int_{\Omega} \MM{w} \cdot \nabla p \diff x
=  \int_{\Omega} \MM{w} \cdot \rho_{0} \begin{pmatrix}
fv\\
b\\
\end{pmatrix} \diff x + \int_{\Omega} \nabla \cdot \MM{w} \, p \diff x,  
\ \ \forall\MM{w}\in  \mathring{\mathbb{V}}_1,
\end{eqnarray}
where we have integrated the pressure gradient term by parts in the third line. 
We then calculate $\eta$ from $\MM{q}$ using \eqref{gi_v} and assemble it over 
the domain. 
If the geostrophic imbalance in the model tends to zero with decreasing $\beta$, 
then the limit is a solution of the SG equations. 
This analysis is performed in section \ref{results_sg}.

\section{Results}\label{results}
In this section, we present the results of the frontogenesis experiments using the 
Eady vertical slice model developed in this study, with the use of the finite element 
code generation library Firedrake \citep{rathgeber2016firedrake}.
The constants used to set up the experiments are shown in section \ref{constants}. 
At the beginning of each experiment, the model is initialised in the way described 
in section \ref{initialisation},
then integrated for 25 days in each experiment. 

The model resolution is given by
\begin{eqnarray}
\Delta x = \frac{2L}{N_{x}},\ \Delta z = \frac{H}{N_{z}}, \label{grid_space}
\end{eqnarray} 
where $N_x$ and $N_z$ are the number of quadrilateral elements in the $x$- 
and $z$-directions, respectively. 
Unless stated otherwise, we use a resolution of $N_{x}$ = 60 and $N_{z}$ = 30, 
which is comparable in terms of DoFs to that 
used in previous work of 
\citet{visram2014framework} and \citet{visram2014asymptotic} 
: $N_{x}$ = 121 and $N_{z}$ = 61.
Note that in our model the effective grid spacings are half the size of 
the lengths given by \eqref{grid_space} as we used a higher-order 
finite element spaces with $k$ = 2, 
as shown in Figure \ref{fig:space-nodes}b and Figure \ref{fig:vertical-nodes}b, 
for all experiments. \citet{nakamura1989nonlinear} used 
a lower resolution of 
$N_{x}$ = 100 and $N_{z}$ = 20. They repeated the experiment with 
twice the horizontal and vertical resolution (results not shown) 
and found very small differences to the low-resolution results.
Therefore we use their result with $N_{x}$ = 100 and $N_{z}$ = 20 
as a comparable result to our control-run result in this section.

For the control run, the time-centring parameter $\alpha$ and the rescaling 
parameter $\beta$ are set to 0.5 and 1, respectively, and a time step of 
$\Delta t$ = 50 s is used based on stability requirements.
In section \ref{results_general}, we investigate the general results of frontogenesis 
from the control run. 
Then the asymptotic convergence of the model to the SG limit is examined in section 
\ref{results_sg}, by repeating the experiment with various $\beta$. 
We also assess the effect of off-centring on the long term performance of the model 
by increasing $\alpha$. 
Finally, we discuss the model results in terms of energy dynamics in section 
\ref{results_energy}.

\subsection{General results of frontogenesis}\label{results_general}
\begin{figure}[p]
  \centering
  \begin{subfigure}{0.45\hsize}
    \centering
    \includegraphics[width=85mm]{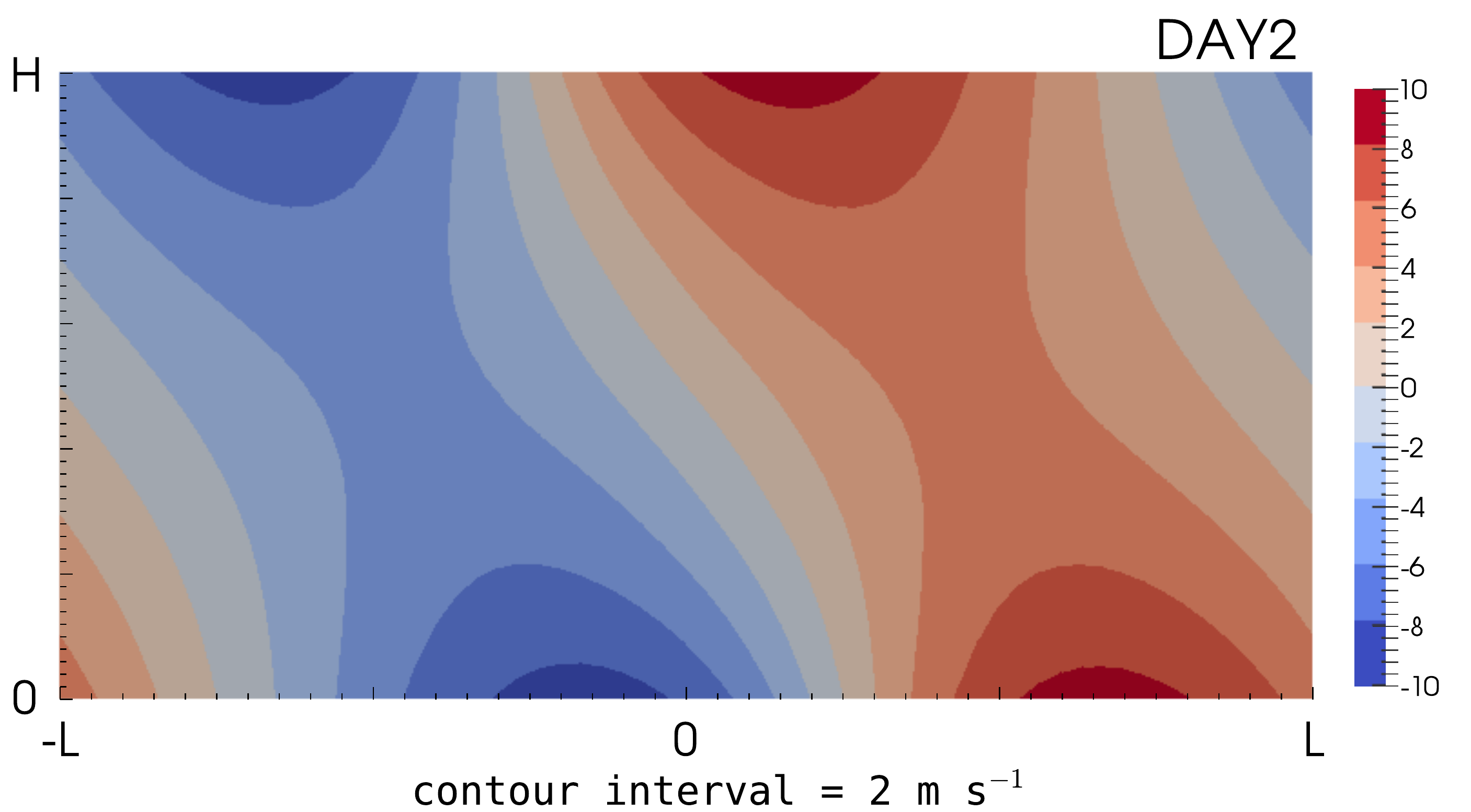}
    \includegraphics[width=85mm]{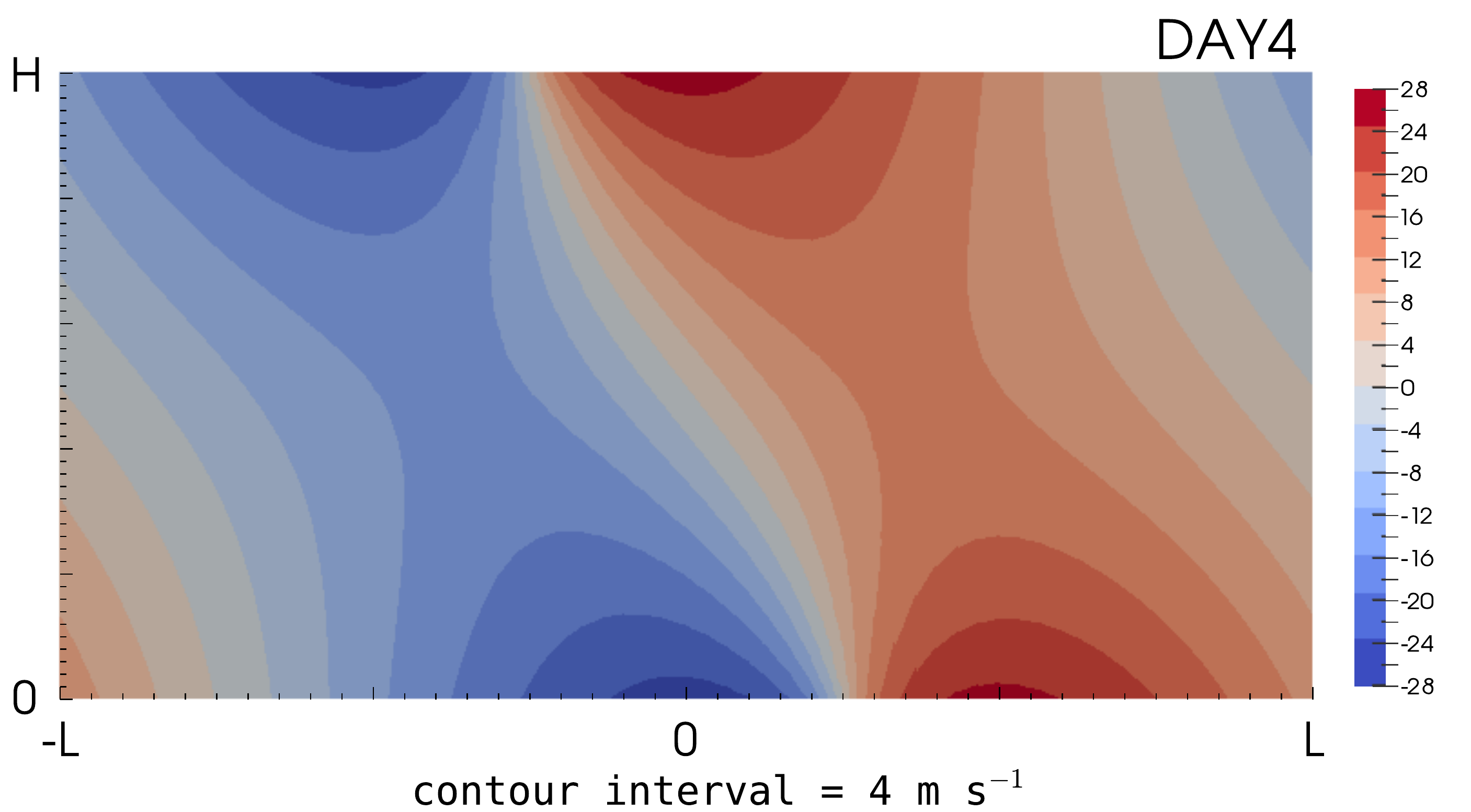}
    \includegraphics[width=85mm]{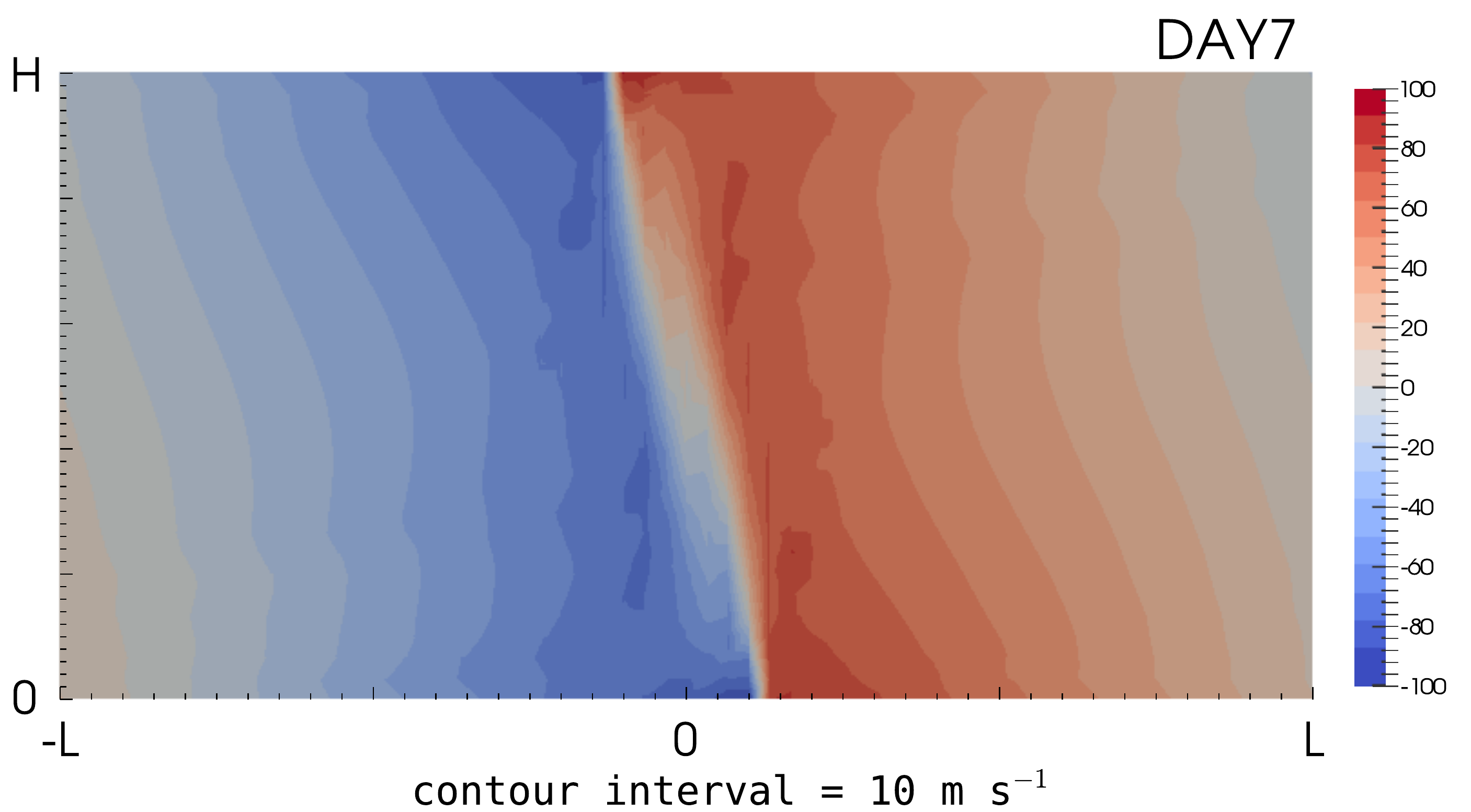}
    \includegraphics[width=85mm]{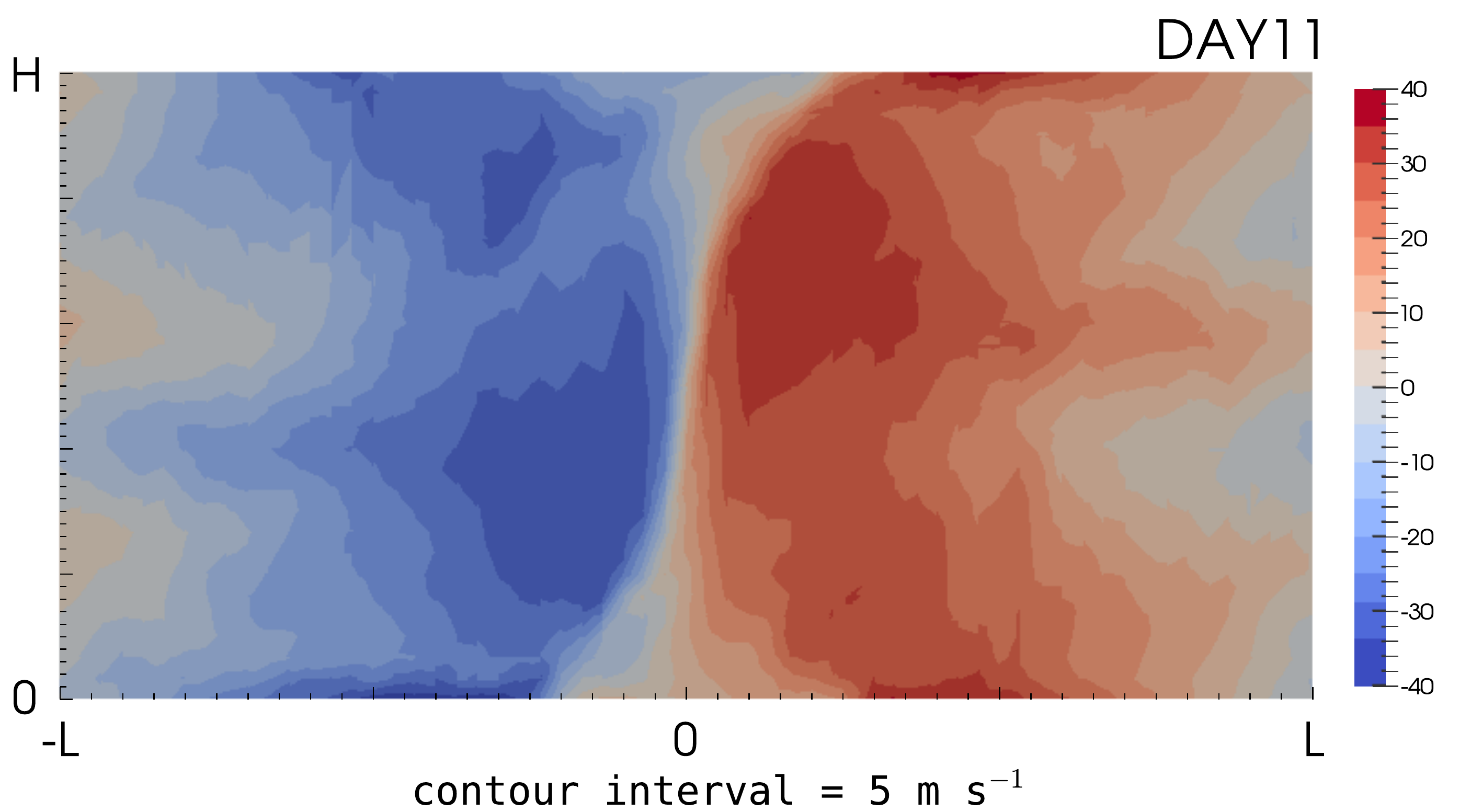}
    \caption{Out-of-slice velocity field. The contour intervals are as specified 
    in each panel.}
    \label{fig:velocity-contours}
  \end{subfigure}
  \hspace{2em}
  \begin{subfigure}{0.45\hsize}
    \centering
    \includegraphics[width=85mm]{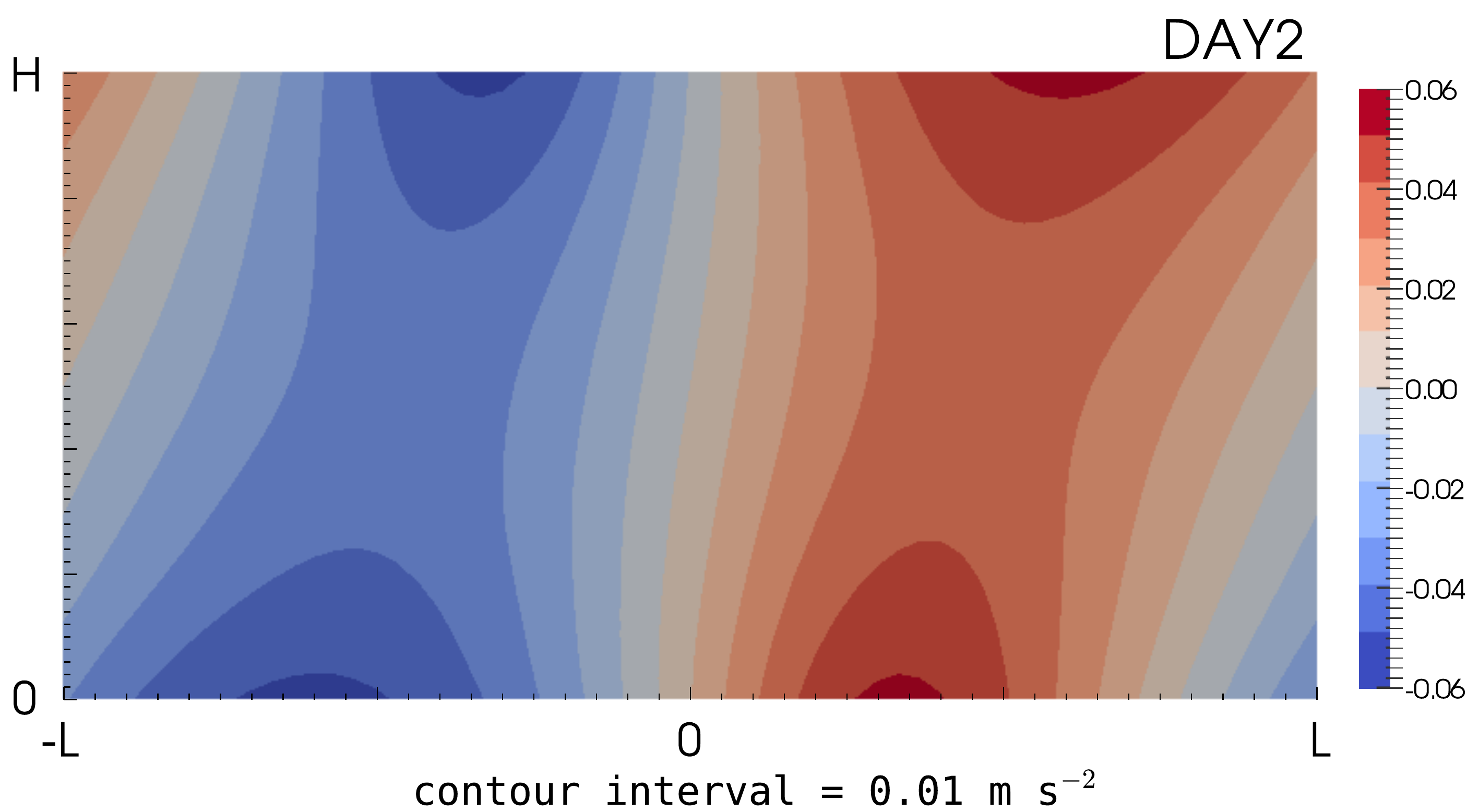}
    \includegraphics[width=85mm]{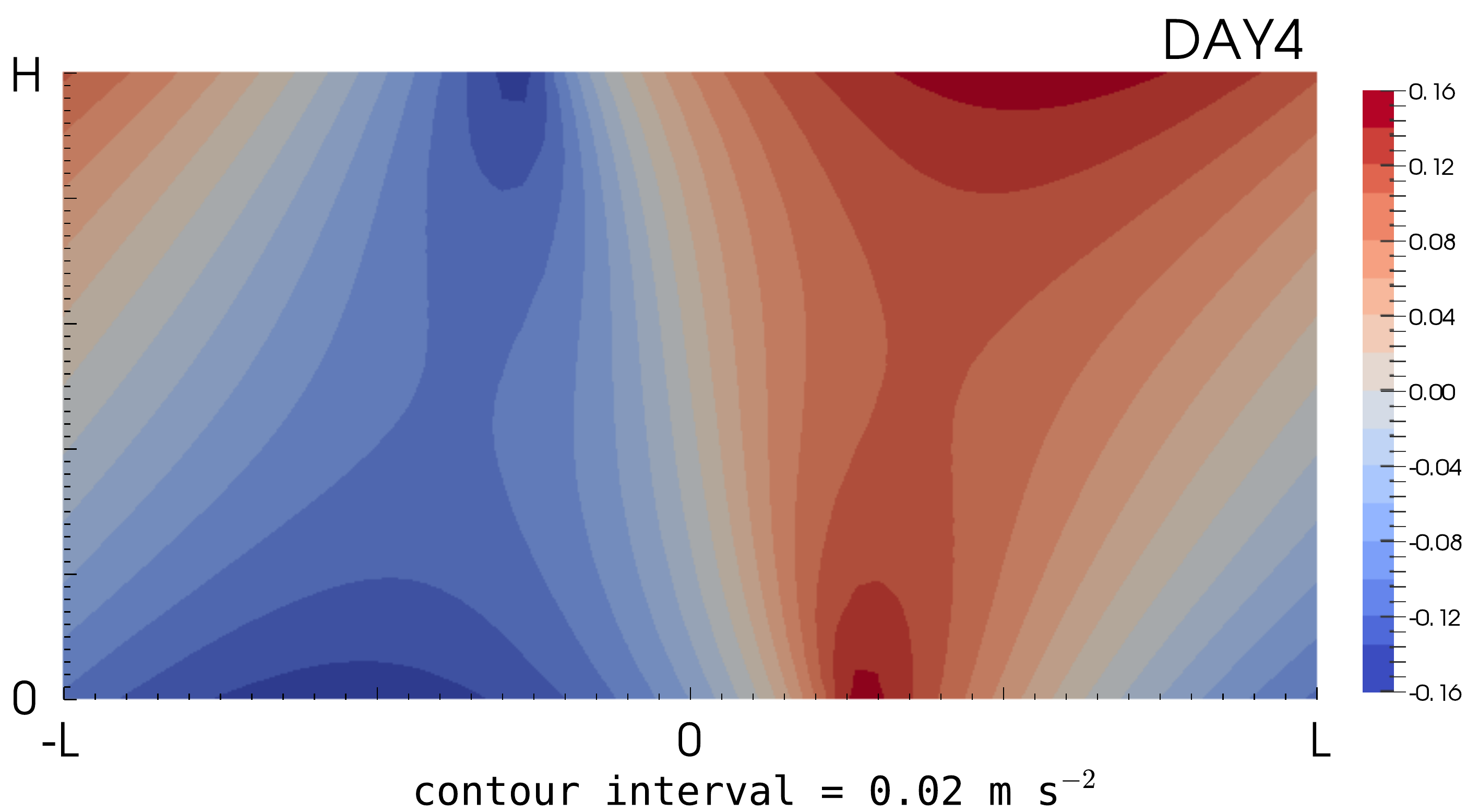}
    \includegraphics[width=85mm]{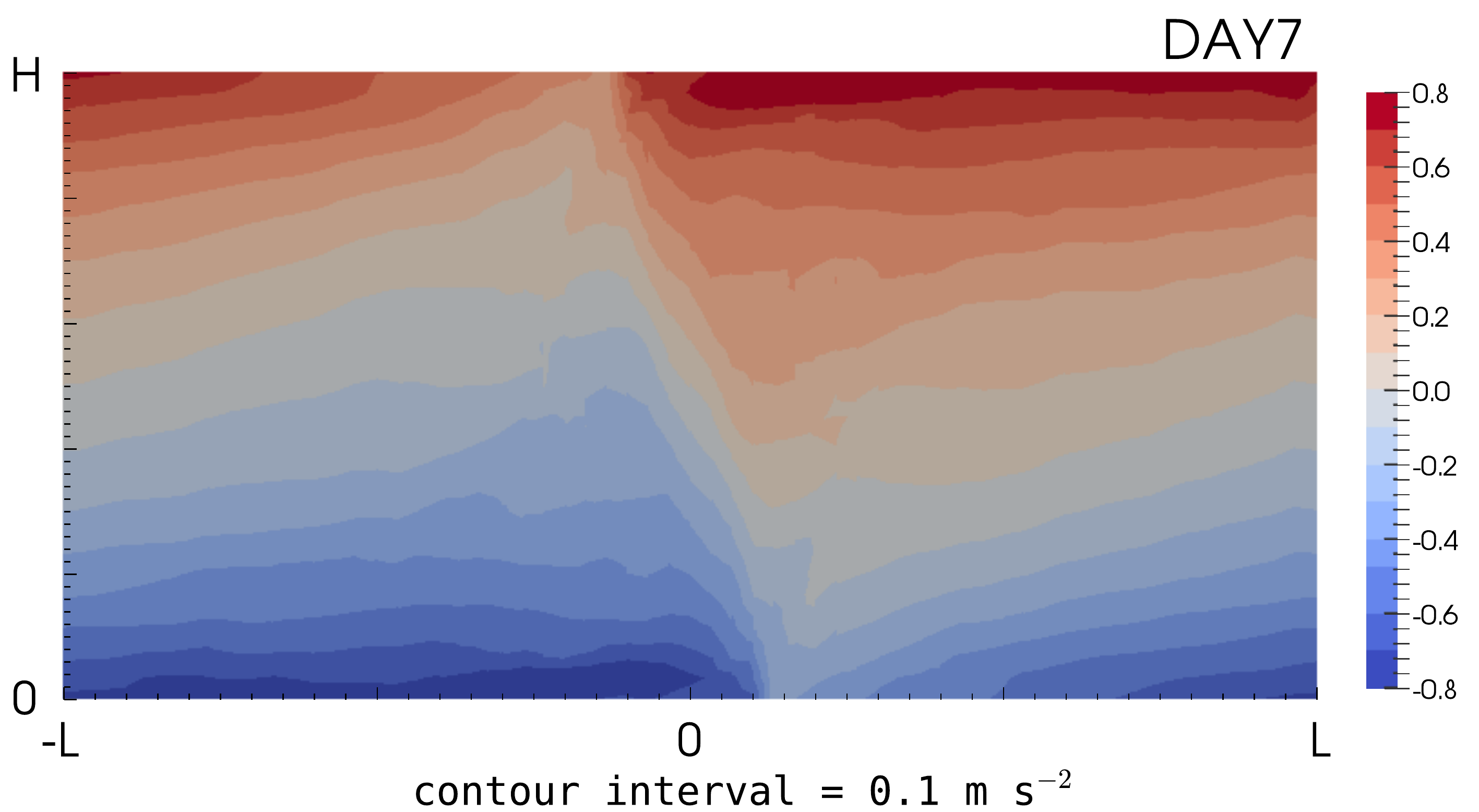}
    \includegraphics[width=85mm]{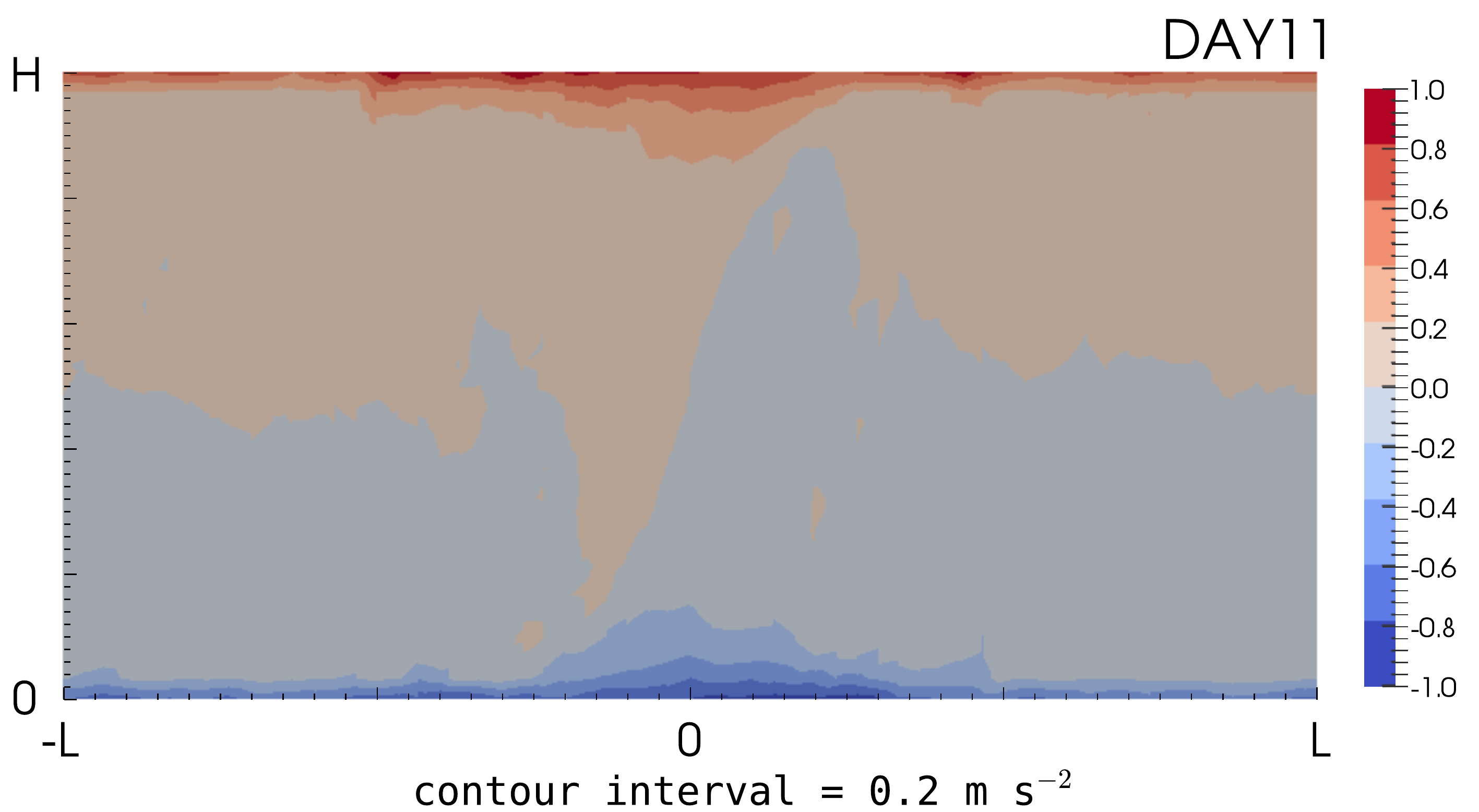}
    \caption{In-slice buoyancy field. The contour intervals are as specified 
    in each panel.}
    \label{fig:buoyancy-contours}
  \end{subfigure}
  \caption{Snapshots of out-of-slice velocity and in-slice bouyancy 
    in the control run at days 2, 4, 7, and 11.}
  \label{fig:slice-contours}
\end{figure}

Figures \ref{fig:velocity-contours} and \ref{fig:buoyancy-contours} show the snapshots 
of the out-of-slice velocity $v$ and buoyancy $b$ fields, respectively, of the control run. 
At day 2, both fields show very similar structures to those from the simulation using 
the linearised equations in \citet{visram2014asymptotic}. 
It suggests that at this early stage the motion is well described by the linearised 
equations. 

The model shows some early signs of front formation at day 4. 
The general shape of the $v$-field is similar to that of day 2. 
However, the gradient in the cyclonic region is now larger than that in the anticyclonic 
region, which indicates the beginning of front formation. 
The maximum gradient in the $v$-field is found at the upper and lower boundaries. 
In the $b$-field, the region of warm air occupies a smaller area at the lower boundary 
than it does at the upper boundary. 
As in the case with $v$, the $b$-field shows the largest gradients near the upper and 
lower boundaries. 

The frontal discontinuity becomes most intense around day 7. 
Strong gradients are found in both $v$ and $b$ fields. 
The frontal zone tilts westward with height in both fields. In the $b$-field, the warm 
region is now lifted off the surface, showing that the front is occluded.

Day 11 corresponds to the first minimum after the initial frontogenesis. 
At this stage the vertical tilt in the $v$-field reverses, which is a sign of energy 
conversion from kinetic back to potential energy. 
In the $b$-field, the discontinuity vanishes and the solution looks almost 
vertically stratified. 

Overall, these results are qualitatively consistent with the early studies 
\citep[e.g.][]{williams1967atmospheric, nakamura1989nonlinear, nakamura1994nonlinear, 
cullen2007modelling, budd2013monge, visram2014framework, visram2014asymptotic}. 
Now, as the out-of-slice velocity is the dominant source of the kinetic energy in the 
Eady problem, we take it as a quantity to compare the strength of the fronts reproduced 
in the models. 

\begin{figure}[htbp]
  \centering
  \includegraphics[width=120mm,angle=0]{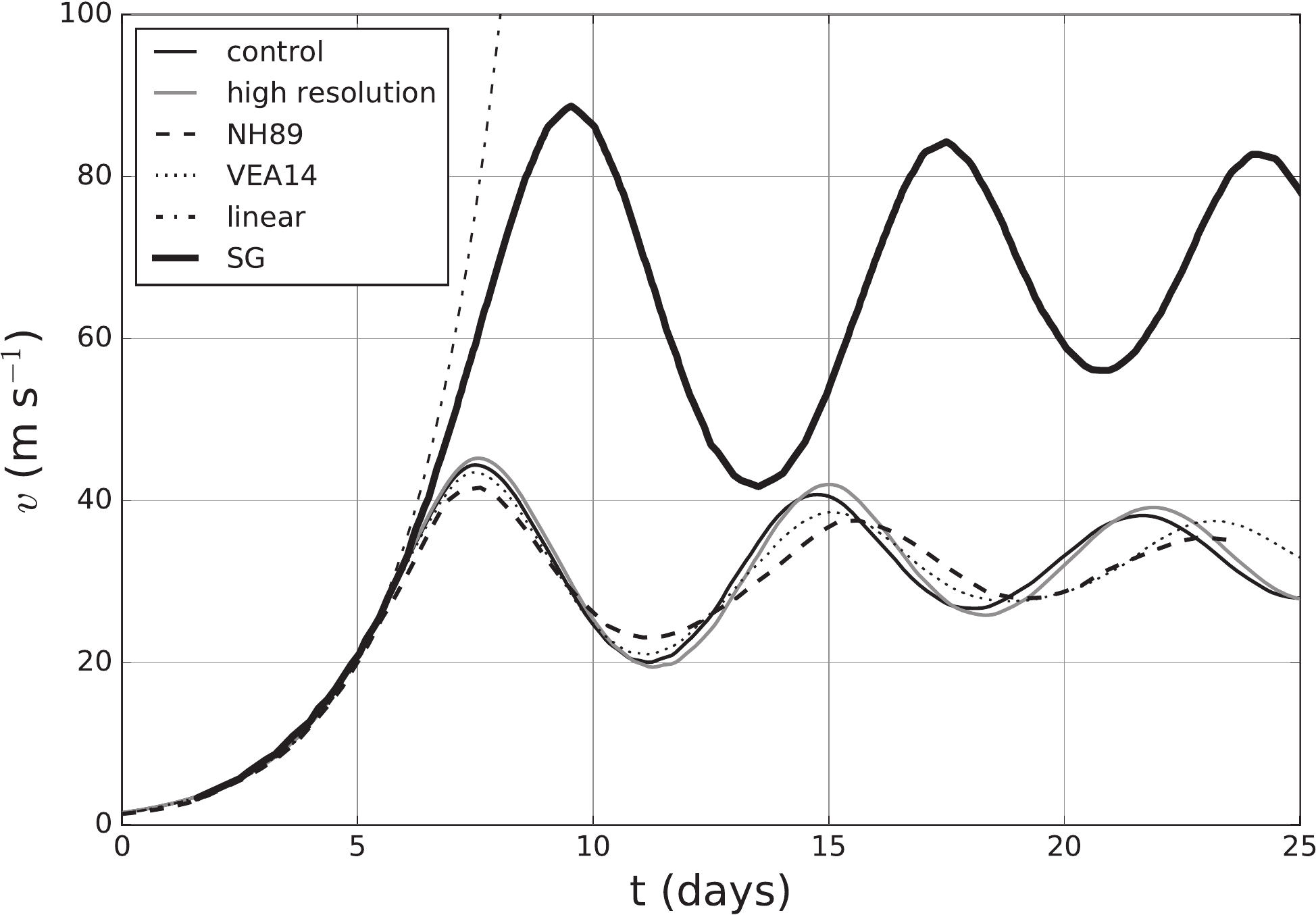}
  \caption{Comparison of the root mean square of the out-of-slice velocity in Eady models. 
  Thin dark solid line shows the result from the control run of this study. 
  Gray line shows the result using two times higher resolution than that of the control run. 
  Dashed and dotted lines show the nonlinear results of 
  \citet{nakamura1989nonlinear} and \citet{visram2014framework}, respectively, 
  and the thick line shows the semi-geostrophic limit solution given by 
  \citet{cullen2007modelling}, based on data from Figure 4 of \citet{visram2014framework}. 
  Dot-dashed line shows the linear result of \citet{visram2014asymptotic} based on data
  from Figure 5.2a of \citet{visram2014asymptotic}.}
  \label{fig:rmsv-velocity-comparison}
\end{figure}

The thin solid curve in Figure \ref{fig:rmsv-velocity-comparison} shows the time 
evolution of the root mean square of $v$ (RMSV) in our model. 
The result shows that the model reproduces several further quasi-periodic lifecycles 
after the first frontogenesis. 
Also shown in Figure~\ref{fig:rmsv-velocity-comparison} in black are the nonlinear results of 
\citet{nakamura1989nonlinear} and \citet{visram2014framework}, the linear result of 
\citet{visram2014asymptotic}, and the SG limit solution from \citet{cullen2007modelling}. 
All results show a good agreement up to around day 5, where the RMSV of the nonlinear 
results grow exponentially following the growth of the linear mode. 
After day 5, at which point the front is close to the grid scale, the nonlinear effects 
become significant and begin to reduce the growth rate. 

For the period of the first frontogenesis,
our result is reasonably close to the result of \citet{visram2014framework}, who applied 
a finite difference method with semi-implicit time-stepping and semi-Lagrangian transport
on the same governing equations as in this study.
Then the two solutions diverge for the subsequent lifecycles. 
Compared to the result of \citet{nakamura1989nonlinear}, who used hydrostatic 
primitive equations with a viscous Eulerian method, both our result and the result of 
\citet{visram2014framework} show larger peak amplitudes of the fronts. 
However, compared to the SG limit solution given by \citet{cullen2007modelling}, 
our result, and the results of \citet{nakamura1989nonlinear} and 
\citet{visram2014framework}, are all much smaller in amplitude.
We believe this to be because \citet{cullen2007modelling} used a Lagrangian
discretisation that resolves fronts even at very coarse resolution, 
and very fine resolution of the front is required to allow this additional transfer of 
potential to balanced kinetic energy. \citet{visram2014framework} showed that the 
Lagrangian conservation properties were badly violated in the Eulerian calculations, 
even at higher resolution, concluding that \cite{cullen2007modelling} is capturing 
the correct solution after the front is formed.

To evaluate the effect of resolution on our model result, 
we repeated the experiment using two times higher resolution than that of the control run: 
($N_x$, $N_z$) = (120, 60). A time step of $\Delta t$ = 25 s is used for the high-resolution run.
The evolution of RMSV in the high-resolution run is shown by the gray curve in Figure 
\ref{fig:rmsv-velocity-comparison}. 
Only a slight increase in the peak amplitudes of RMSV is found in the high-resolution 
run compared to that of the control run. 
It indicates that, due to the rapid formation of the frontal discontinuity, the front reaches 
the grid scale very quickly even when double the resolution is used. 
As a result, the use of the high resolution can only slightly delay the collapse of the fronts,
and thus makes very little contribution to filling the gap between the RMSV values of the model 
and the SG limit. 
This result suggests that we would need resolution 
of several orders of magnitude greater than presently used
to reach the RMSV of the SG limit.

\begin{figure}[htbp]
  \centering
  \begin{subfigure}{0.8\hsize}
    \centering
    \includegraphics[width=120mm]{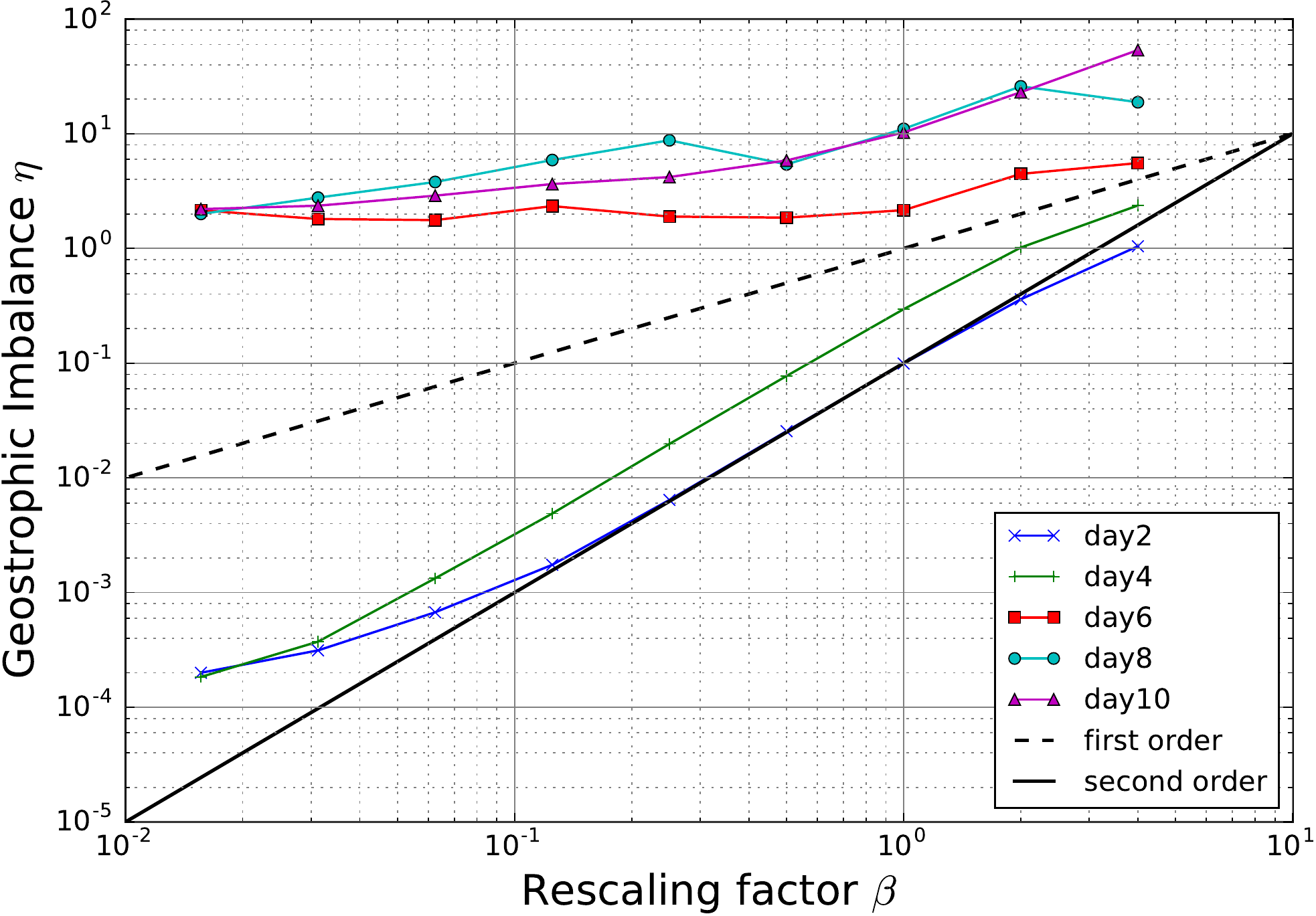}
    \caption{$\alpha = 0.5$}
    \label{fig:geostrophic-imbalance-a}
  \end{subfigure}
  
  \vspace{1em}
  \begin{subfigure}{0.8\hsize}
    \centering
    \includegraphics[width=120mm]{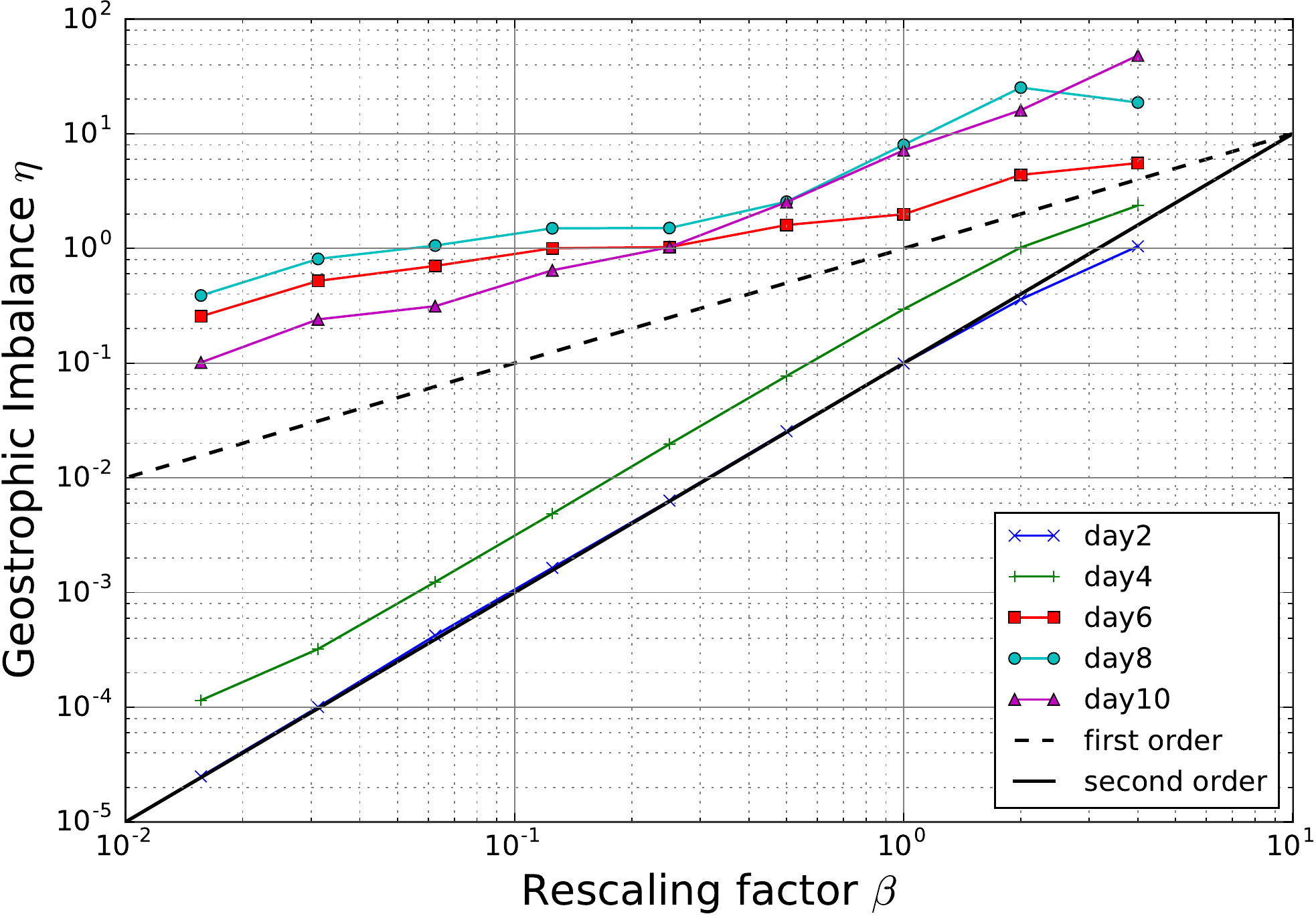}
    \caption{$\alpha = 0.55$}
    \label{fig:geostrophic-imbalance-b}
  \end{subfigure}
  \caption{Comparison of the geostrophic imbalance in the results of the rescaling 
  tests using (a) $\alpha$ = 0.5, and (b) $\alpha$ = 0.55. 
  Black dashed and solid lines correspond to the first- and second-order 
  convergence rates, respectively. 
  Colored lines show the variations of the geostrophic imbalance with rescaling 
  parameter $\beta$ at day 2, 4, 6, 8  and 10.}
  \label{fig:geostrophic-imbalance}
\end{figure}

\begin{figure}[htbp]
  \centering
  \begin{subfigure}{0.8\hsize}
    \centering
    \includegraphics[width=120mm]{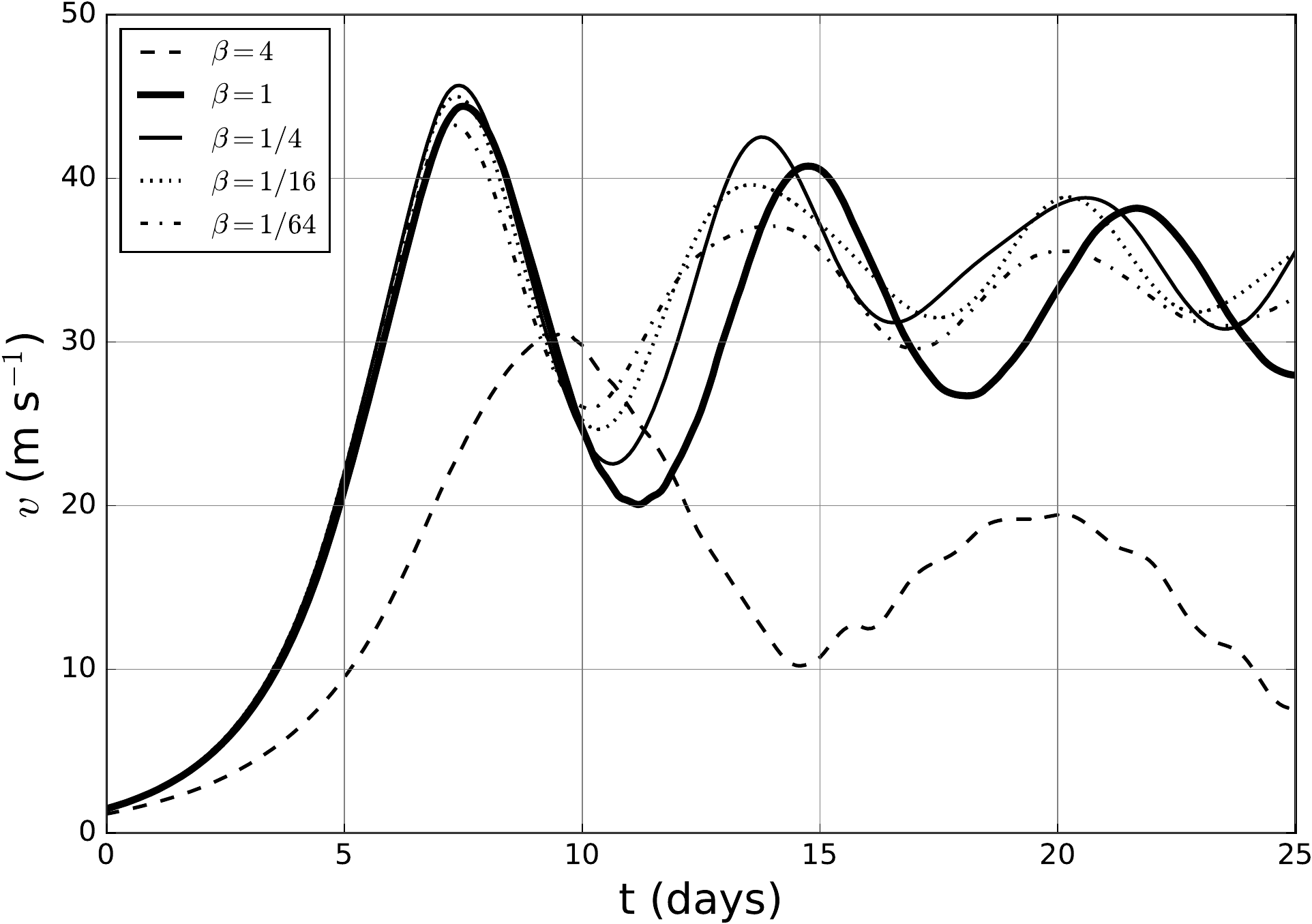}
    \caption{$\alpha = 0.5$}
    \label{fig:rmsv-velocity-rescaling-a}
  \end{subfigure}
  
  \vspace{1em}
  \begin{subfigure}{0.8\hsize}
    \centering
    \includegraphics[width=120mm]{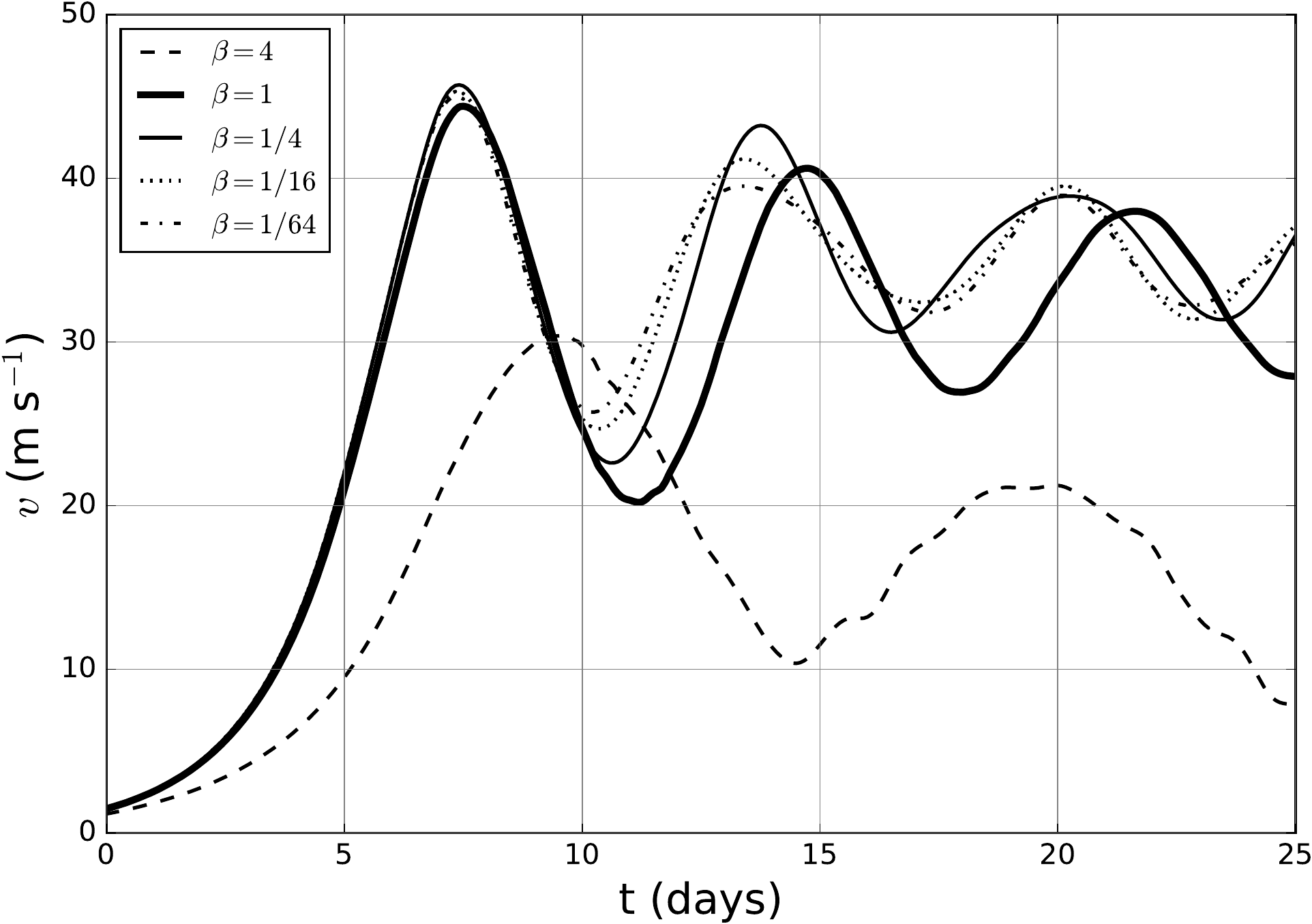}
    \caption{$\alpha = 0.55$}
    \label{fig:rmsv-velocity-rescaling-b}
  \end{subfigure}
  \caption{Comparison of the root mean square of the out-of-slice velocity in the results 
  of the rescaling tests using (a) $\alpha$ = 0.5, and (b) $\alpha$ = 0.55. 
  Thick lines show the results with $\beta$ = 1. 
  The other lines show the results with different values of $\beta$ as shown in the legend.}
  \label{fig:rmsv-velocity-rescaling}
\end{figure}

\subsection{Asymptotic convergence to the SG solution}\label{results_sg}
In this section, the validation test of the asymptotic convergence outlined in section 
\ref{settings_sg} is performed. First, the frontogenesis experiment is repeated 
using eight different values of the rescaling parameter: 
$\beta$ = 2$^{2}$, 2, 2$^{-1}$, 2$^{-2}$, 2$^{-3}$, 2$^{-4}$, 2$^{-5}$, and 2$^{-6}$.
The time step $\Delta t$ is set to 50 s for the experiments 
using 2$^{-2}$ $\leq$ $\beta$ $\leq$ 2$^{2}$, 
25 s for $\beta$ = 2$^{-3}$ and 2$^{-4}$, 
and 12.5 s for $\beta$ = 2$^{-5}$ and 2$^{-6}$.
The other settings are the same as the control run.
We calculated the geostrophic imbalance $\eta$ defined by \eqref{gi} in each experiment. 
We then plotted it over the period of the initial frontogenesis, together with  $\eta$ in the 
control run where $\beta$ is unity. 

Figure \ref{fig:geostrophic-imbalance-a} shows the variation of the geostrophic imbalance 
$\eta$ with rescaling parameter $\beta$ alongside the theoretical first- and second-order 
convergence rates. 
The convergence rate starts at around the second-order for $\beta \geq 2^{-3}$ and 
the first-order for $\beta < 2^{-3}$ as shown by the slope at day 2. 
It improves to the second-order for $\beta \geq 2^{-5}$ at day 4. 
However, it doesn't converge at all at day 6, at which point a strong discontinuity is formed 
in the model. 
After the peak of the initial frontogenesis at around day 7, the convergence rate recovers 
a little but stays at less than first-order at days 8 and 10. 

For the results in Figure \ref{fig:geostrophic-imbalance-a}, the time-centring parameter 
$\alpha$ is set to 0.5 as that is in the control run. 
To test the effect of off-centring, the rescaling test was repeated using $\alpha$ = 0.55. 
This result is shown in Figure \ref{fig:geostrophic-imbalance-b}. 
It shows that increasing the implicitness of the solution gives more balanced solutions. 
In particular, a reduction of the imbalance is found throughout the initial frontogenesis 
for $\beta < 2^{-3}$.
As a result, the overall second-order convergence is achieved at day 2, and the 
first-order convergence is recovered at day 10 in Figure \ref{fig:geostrophic-imbalance-b}. 
This result is comparable to the result of \citet{visram2014framework} (see their Figure 2), 
where $\alpha =$ 0.55 was used,
indicating that the compatible finite element method is performing as well as a finite 
difference method for this test problem. Note that \citet{visram2014framework} used the 
range of $2^{-3} \leq \beta \leq 2^{2}$, whereas we show the convergence of 
geostrophic imbalance for smaller $\beta$ as well. 
The result is also consistent with the results reported by \citet{cullen2007modelling} 
with compressible equations. 

Figure \ref{fig:rmsv-velocity-rescaling-a} and \ref{fig:rmsv-velocity-rescaling-b} 
show the evolutions of RMSV in each experiment corresponding to Figure 
\ref{fig:geostrophic-imbalance-a} and \ref{fig:geostrophic-imbalance-b}, respectively. 
In Figure \ref{fig:rmsv-velocity-rescaling-b}, there are some increase in the peak 
amplitude of fronts compared to that in Figure
\ref{fig:rmsv-velocity-rescaling-a}, especially in the second peak of the results 
for $\beta \leq 2^{-2}$.
It indicates that damping out some of the unbalanced motion with off-centring 
improves the predictability of quasi-periodic lifecycles.
Compared to the off-centred results of \citet{visram2014framework} (see their Figure 4),
which very quickly began to diverge as they decreased the Rossby number,
our model shows 
good predictability throughout the range of $\beta$.
However, in both cases of Figure \ref{fig:rmsv-velocity-rescaling}, 
decreasing $\beta$ does not make a big difference in the peak amplitude of fronts, 
thereby leaving the large gap 
between the model result and the SG limit solution on the strength of the fronts 
regardless of the value of $\beta$.

\subsection{Energy dynamics}\label{results_energy}
In sections \ref{results_general} and \ref{results_sg}, we showed that our model 
reproduced results which are consistent with the early model studies based on finite
difference methods, and showed that the model solutions converge to 
a solution in geostrophic balance when we decrease the Rossby number. 
In this section, we will again look into our result of the control run with focus on 
the energy dynamics.

\begin{figure}[htbp]
  \centering
  \includegraphics[width=120mm]{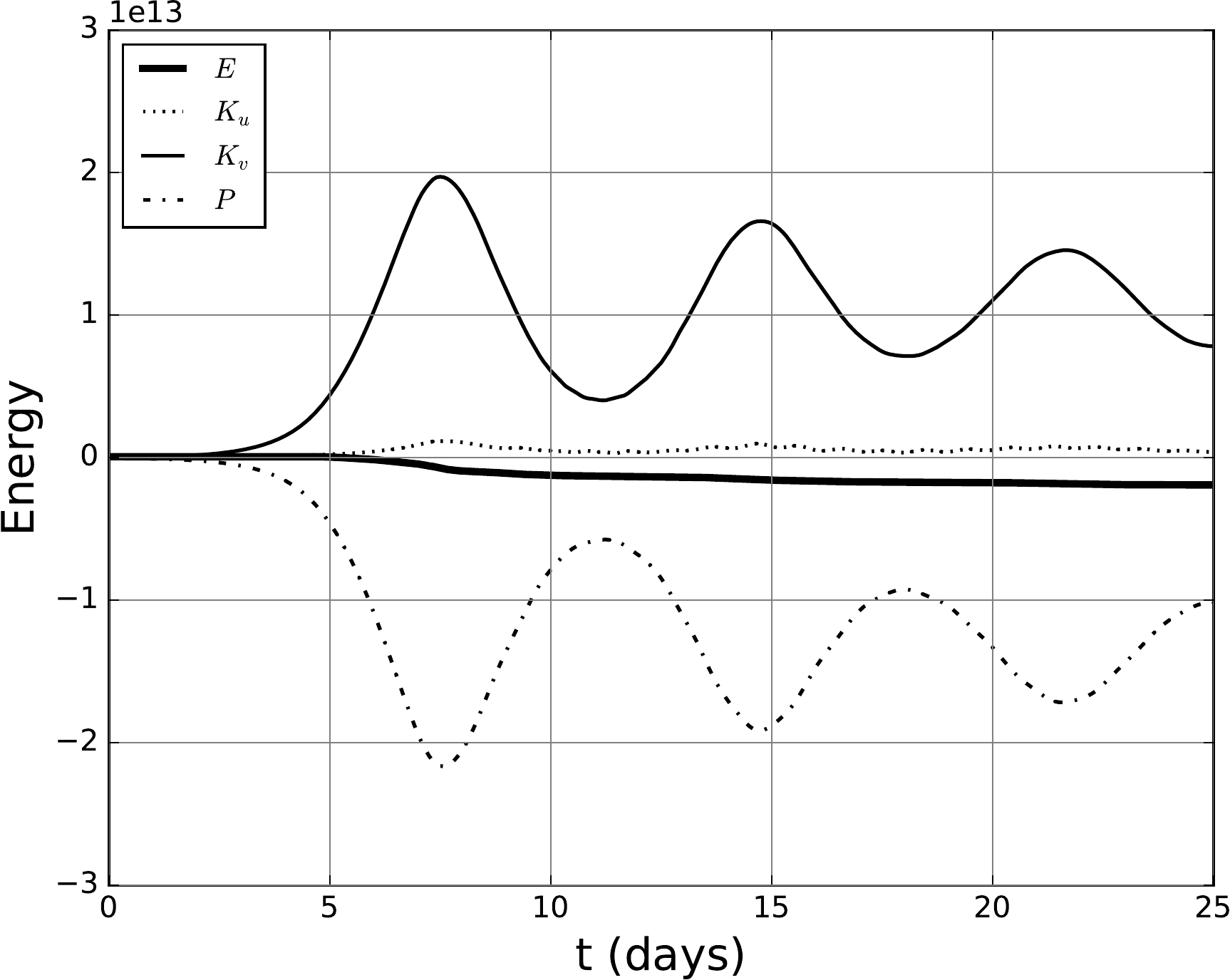}
  \caption{Time evolution of energy of the control run. 
  Thick line represents the evolution of total energy. 
  Dotted and solid lines represent the evolutions of in-slice and out-of-slice 
  components of the kinetic energy. 
  Dot-dashed line represents the evolution of potential energy.}
  \label{fig:time-energy-control}
\end{figure}

Figure \ref{fig:time-energy-control} shows the time evolution of the total energy $E$, 
the kinetic energy $K_{u}$ and $K_{v}$, and the potential energy $P$, which are 
defined by the equations \eqref{total_energy} to \eqref{potential}. 
The kinetic energy $K_{v}$ reaches the maximum amplitude at around day 7, then 
reduces to the first minimum at around day 11 followed by smaller amplitude 
lifecycles, just as RMSV does in Figure \ref{fig:rmsv-velocity-comparison}.
The time evolution of potential energy $P$ shows the same behaviour with 
opposite sign, which demonstrates the exchange from potential to kinetic energy 
over several lifecycles. 
The amplitude of the kinetic energy $K_u$ is very small compared to that of $K_v$ 
and $P$ throughout the experiment. 
As a result, the total energy $E$ can be interpreted as the difference of the 
amplitude of $K_{v}$ from that of $P$, which shows a gradual decrease with time.

\begin{figure}[htbp]
  \centering
  \includegraphics[width=120mm]{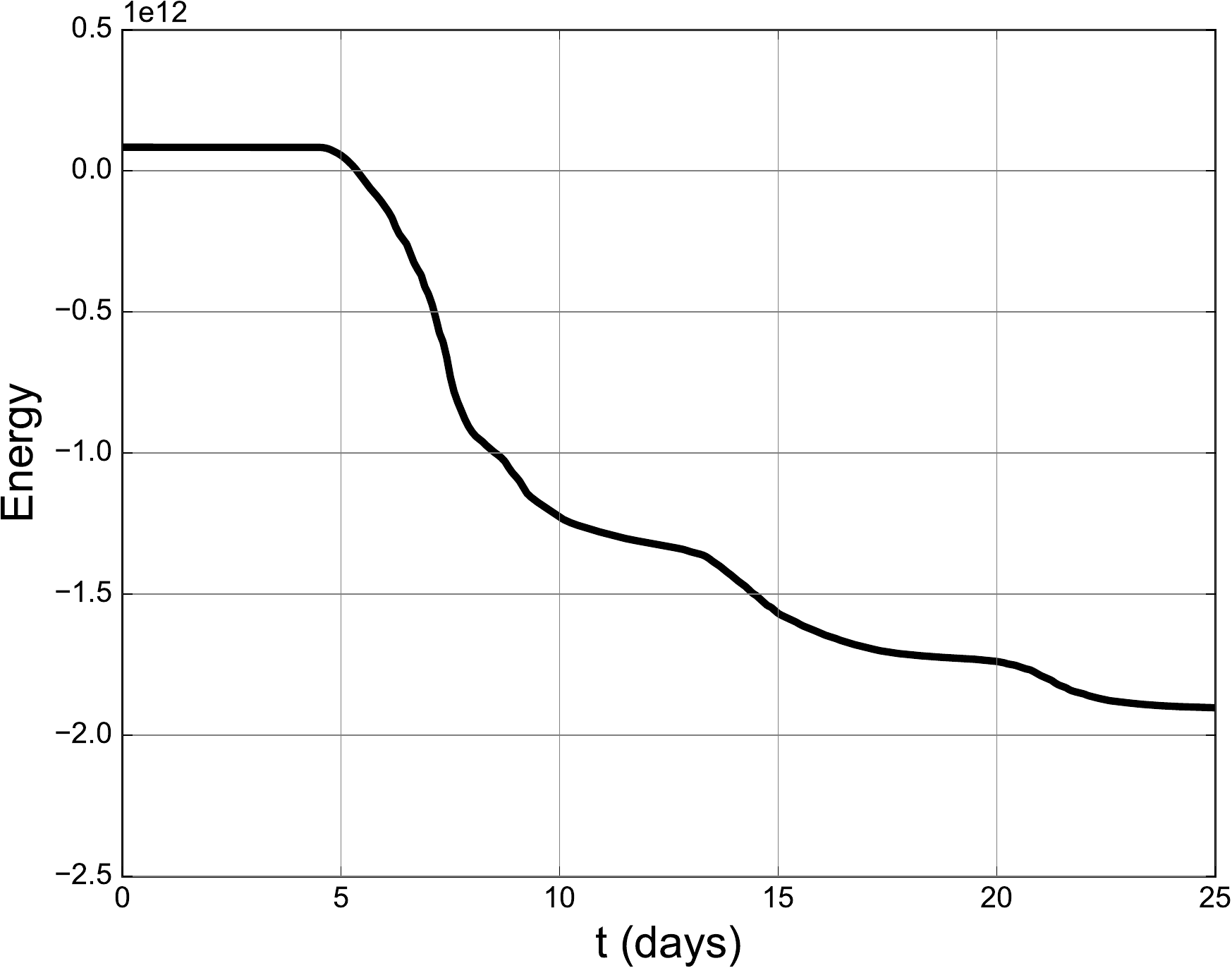}
  \caption{Enlarged view of the evolution of total energy shown 
  in Figure \ref{fig:time-energy-control}.
  The vertical scale is one order of magnitude less than 
  that of Figure \ref{fig:time-energy-control}.}
  \label{fig:time-total-energy}
\end{figure}

Figure \ref{fig:time-total-energy} provides an enlarged view of the time evolution 
of the total energy. 
Note that the thick lines in Figure \ref{fig:time-energy-control} and 
Figure \ref{fig:time-total-energy} show the same evolution of the total energy 
of the control run, and the vertical scale of Figure \ref{fig:time-total-energy} is 
one order of magnitude less than that of Figure \ref{fig:time-energy-control}.
The total energy stays constant up until day 5, then starts decreasing. 
It becomes quasi-constant from around day 10 to day 13, then decreases again. 
By comparing this with the lifecycles of fronts shown as the evolution of RMSV 
in Figure \ref{fig:rmsv-velocity-comparison}, it appears that the model starts 
losing energy every time the discontinuity reaches the grid scale. 
In addition, the reduction in the total energy starts at almost the same time as 
the RMSV of the control run diverges from the SG solution. 
Therefore we consider that the lack of resolution is a significant cause of the 
loss of energy, and it also stops the growth of RMSV too early. 

\begin{figure}[htbp]
  \centering
  \includegraphics[width=120mm]{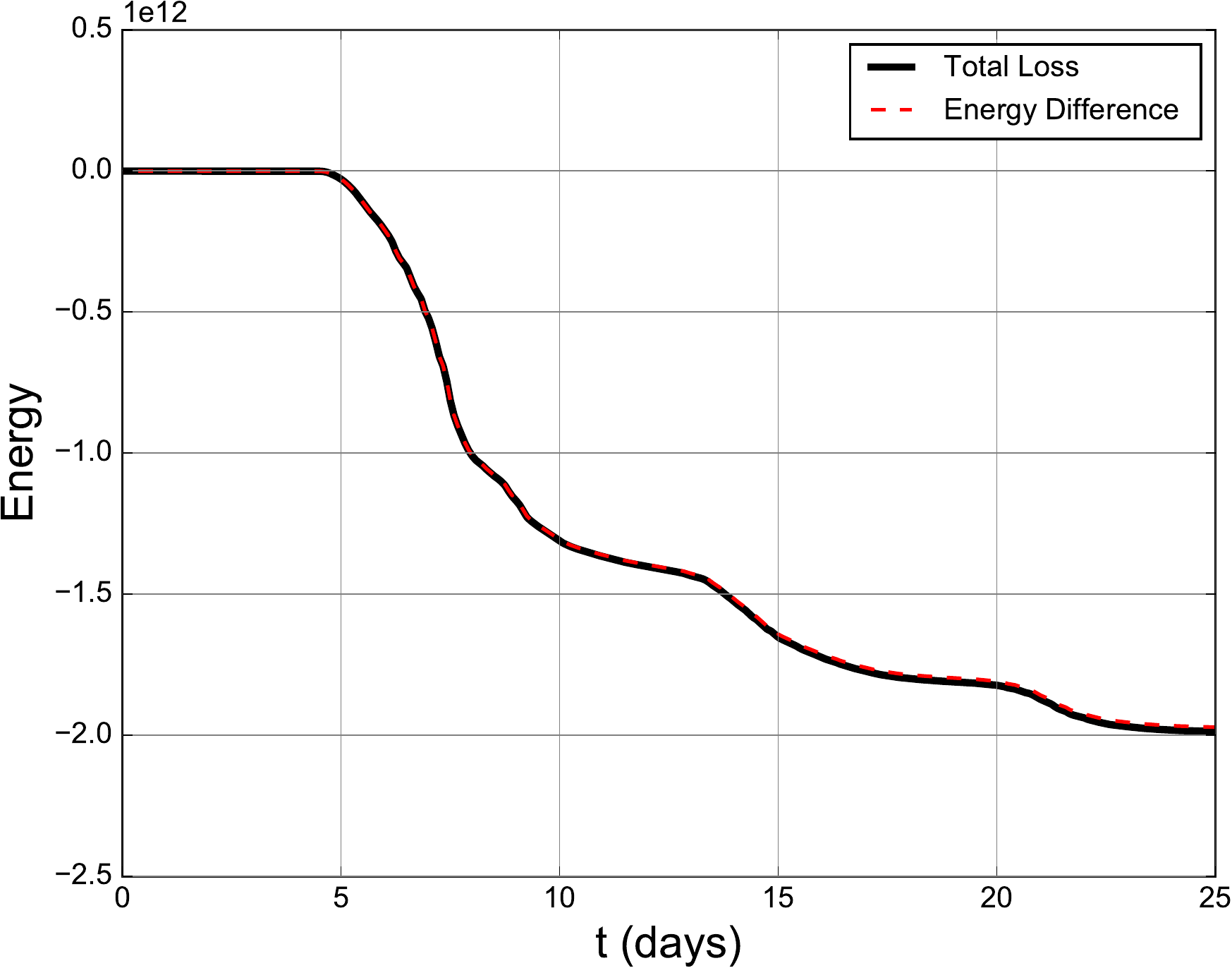}
  \caption{Result of the energy analysis using the dummy velocity. 
  Thick line represents the loss of total energy
  from the initial state of the control run. Dashed
  line shows the accumulated difference between the kinetic
  energy $K_v$ and the dummy kinetic energy $K_{v_d}$.}
  \label{fig:energy-dummy}
\end{figure}

Now, to estimate the potential loss in the kinetic energy $K_{v}$ caused by our 
advection scheme for $v$, we perform a test considering the dummy velocity 
$v_d$ which obeys the following advection-only equation,
\begin{eqnarray}
\frac{\partial v_d}{\partial t} +  \MM{u} \cdot \nabla v_d &=& 0, \label{dummyv}
\end{eqnarray}
and the dummy kinetic energy $K_{v_d}$ defined by
\begin{eqnarray}
K_{v_d} &=& \rho_{0} \int_{\Omega} \frac{1}{2} v_d^{2} \ \diff x.\label{kinetic_vd}
\end{eqnarray}
In this test, we solve the equation \eqref{dummyv} in parallel with the governing 
equations \eqref{ueq} to \eqref{peq}. The same experimental settings 
including the constants, initial and boundary conditions and the resolution 
as in the control run are used in this test. 
Here we apply the same advection scheme used for $v$ to $v_d$ as
\begin{eqnarray}
\int_{\Omega}\phi\frac{\partial v_d}{\partial t} \diff x 
-  \int_\Omega \nabla \cdot (\phi \MM{u}) v_d \diff x 
+ \int_\Gamma \jump{\phi \MM{u}} \tilde{v_d} \diff S 
 = 0, \ \ \forall\phi\in  \mathbb{V}_2,  \label{dummyv_dis}
\end{eqnarray}  
and the same semi-implicit time-stepping method described in section 
\ref{time_discretisation} to \eqref{dummyv_dis}. 
After every time step, we calculate the difference between $K_v$ and 
$K_{v_d}$ as
\begin{eqnarray}
\epsilon =  \rho_{0} \int_{\Omega} \frac{1}{2} \{(v_d^{n+1})^{2} - (v^{n})^{2}\} 
\ \diff x. 
\end{eqnarray} 
Then we replace $v_{d}^{n+1}$ with $v^{n+1}$ and repeat the time integration. 
By accumulating the energy difference $\epsilon$ every time step, we can 
estimate the potential loss of kinetic energy caused by the discretisation of 
the advection term in the $v$ equation \eqref{veq}. 
This result is shown by the dashed line in Figure \ref{fig:energy-dummy}. 
Also shown in Figure \ref{fig:energy-dummy} as the thick line is the loss 
of total energy in the control run from the initial state. 
The two energy values are almost identical to each other, showing that
almost all the energy loss in the model is caused in the $v$ advection. 

\begin{figure}[htbp]
  \centering
  \includegraphics[width=140mm,angle=0]{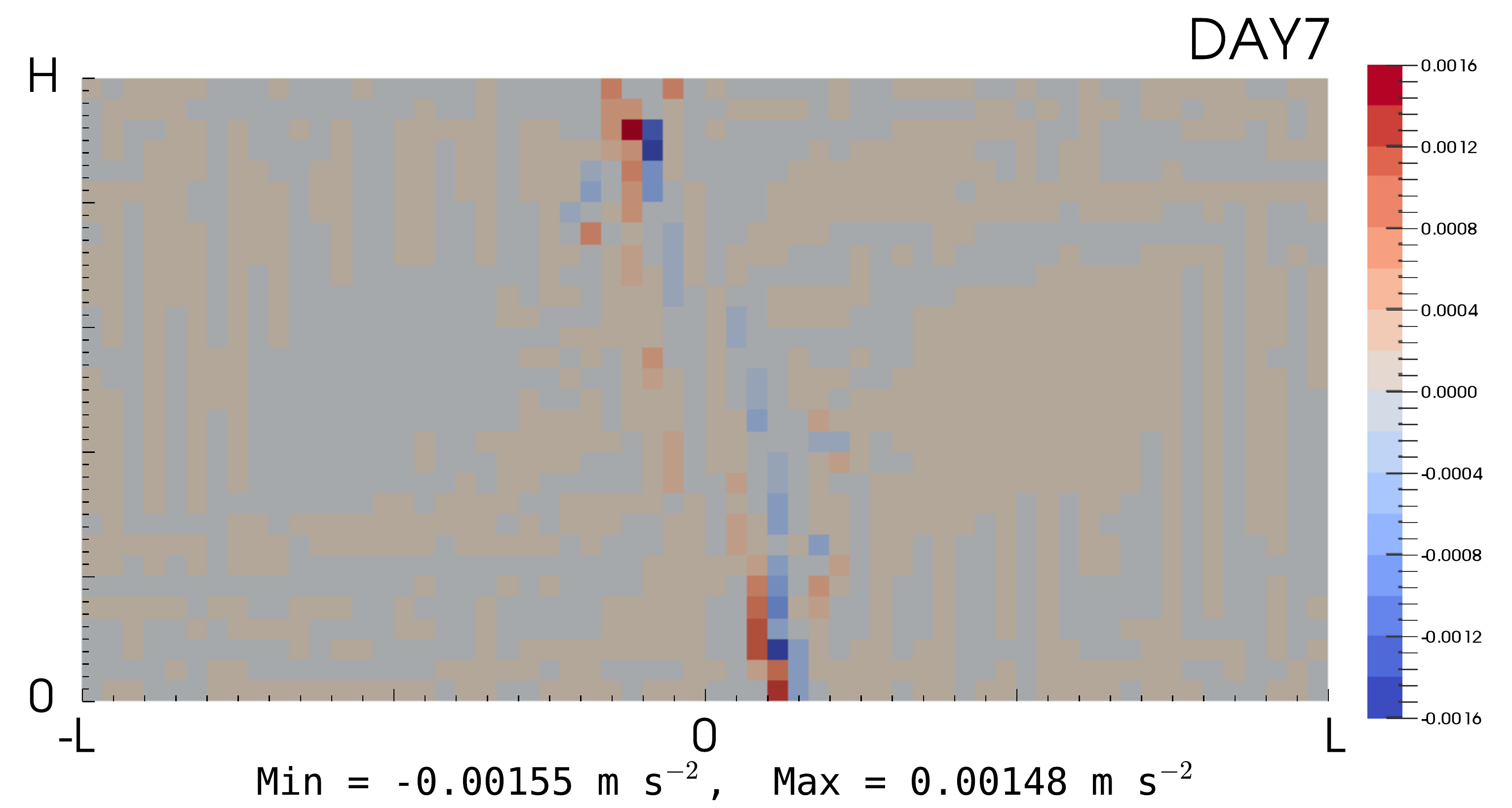}
  \caption{Residual in the out-of-slice velocity field calculated at day 7 of the control run.}
  \label{fig:velocity-slice-residual}
\end{figure}

\begin{figure}[htbp]
  \centering
  \includegraphics[width=120mm,angle=0]{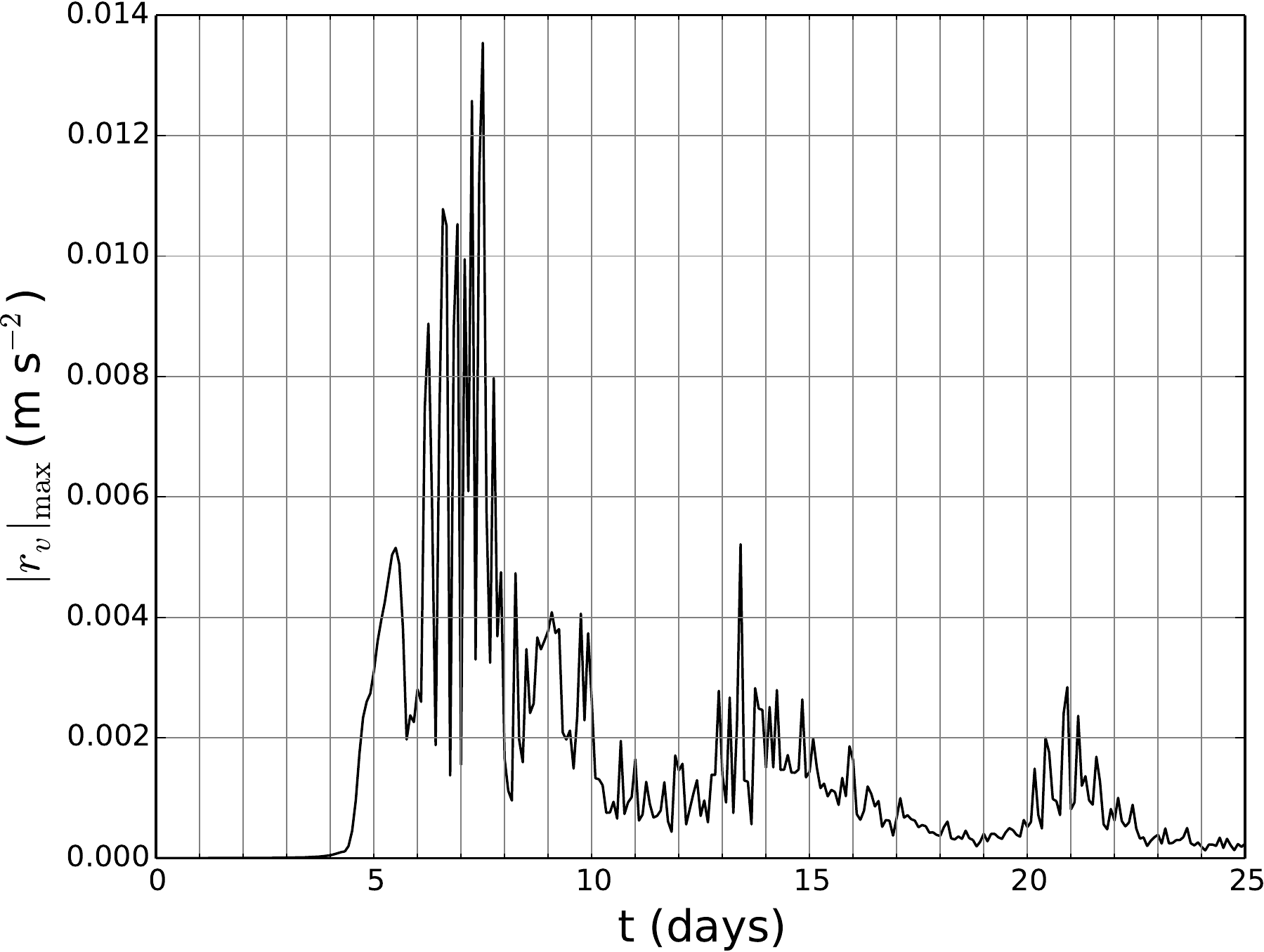}
  \caption{Time evolution of the maximum amplitude of the residual in the out-of-slice 
  velocity field of the control run.}
  \label{fig:time-residual}
\end{figure}

To investigate the cause of the energy loss in the $v$ advection, 
we calculated the residual in the out-of-slice velocity field, which is defined 
as the difference between LHS and RHS of the equation \eqref{veq},
\begin{eqnarray}
r_{v} = \frac{\partial v}{\partial t} +  \MM{u} \cdot \nabla v + f\MM{u} \cdot \hat{\MM{x}} 
+ \frac{\del \bar{b}} {\del y}\left(z-\frac{H}{2}\right).
\end{eqnarray}
With $\phi \in \mathbb{V}_2$, we calculated $r_{v}$ by solving
\begin{eqnarray}
\int_{\Omega} \phi r_{v}^{n+1} \diff x = 
\int_{\Omega} \phi \frac{v^{n+1}-v^{n}}{\Delta t} \diff x 
+ \int_{\Omega} \phi \MM{u}^{n+\frac{1}{2}} \cdot \nabla v^{n+\frac{1}{2}} \diff x 
+ \int_{\Omega} \phi f\MM{u}^{n+\frac{1}{2}} 
\cdot \hat{\MM{x}} \diff x 
+ \int_{\Omega} \phi \frac{\del \bar{b}} {\del y}\left(z-\frac{H}{2}\right) \diff x, 
\ \ \forall\MM{\phi}\in  \mathbb{V}_2. \label{vtres}
\end{eqnarray}
Figure \ref{fig:velocity-slice-residual} shows the residual at day 7, which is 
when the front reaches the first peak. 
Compared with the $v$-field at day 7 in Figure \ref{fig:velocity-contours}, we can see 
a large increase in the residual occurring along the frontal discontinuity. 
In particular, the maximum amplitude of the residual is found near the upper and lower 
boundaries, where the discontinuity is most intense.
Figure \ref{fig:time-residual} shows the time evolution of the maximum amplitude of $r_v$.
It shows the biggest peak during the first lifecycle followed by small peaks during the 
second and third lifecycles.
These results indicate that our advection scheme for $v$ does not converge well enough 
at the fronts due to the strong discontinuity. 
This then leads to the dissipation of the the kinetic energy $K_{v}$ in the model 
every time the discontinuity reaches the grid scale.

\begin{figure}[htbp]
  \centering
  \includegraphics[width=120mm]{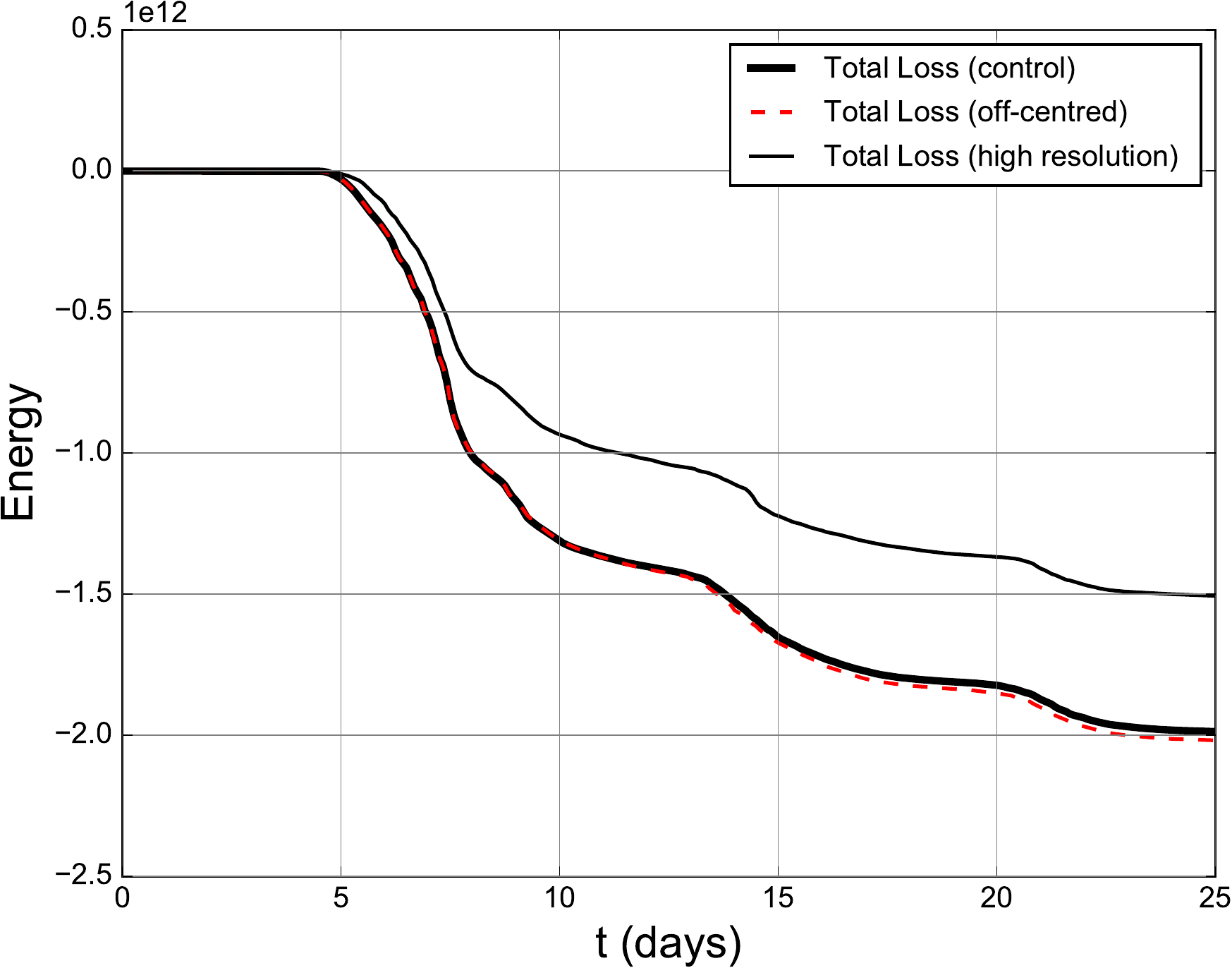}
  \caption{Comparison of the loss of total energy from the initial state
  in experiments with different settings. 
  Thick line represents the same loss of total energy of the control 
  run as in Figure \ref{fig:energy-dummy}.  
  Dashed line shows the loss when using $\alpha$ = 0.55.
  Thin line represents the loss in the high-resolution run.}
  \label{fig:compare-total-energy}
 \end{figure}

Finally, Figure \ref{fig:compare-total-energy} provides a comparison of the loss 
of total energy in the experiments performed in the previous sections. 
Note that with all three results shown in Figure \ref{fig:compare-total-energy} 
the rescaling factor $\beta$ is unity. 
The thick line shows the same loss of total energy from the initial state in the control run 
as in Figure \ref{fig:energy-dummy}.
The dashed line shows the loss of the total energy in the 
experiment with off-centring, which is performed in section \ref{results_sg}. 
It is shown that the use of off-centring has an insignificant effect 
on the loss of energy. 
This result, together with the results in Figures \ref{fig:geostrophic-imbalance-b} 
and \ref{fig:rmsv-velocity-rescaling-b}, suggests that off-centring does not have 
a big impact on the large scale dynamics but the unbalanced motion. 
The thin solid curve in Figure \ref{fig:compare-total-energy} represents 
the loss of total energy in the high-resolution run, which is performed 
in section \ref{results_general}.
There is a clear improvement in energy conservation with the use of the high resolution;
the total loss at day 25 is about 25 \% less than that of the control run. 
It supports our assumption that the lack of resolution is a significant cause 
of the loss of energy.
However, as shown in Figure \ref{fig:rmsv-velocity-comparison},
the high-resolution run gives only a slight increase in the peak amplitudes of RMSV 
despite the improvement in energy conservation.
Therefore we have concluded that the energy loss in the model does not 
account for the large gap between the model result and the SG limit solution
on RMSV.

\section{Conclusion}\label{conclusion}
A new vertical slice model of nonlinear Eady waves was developed using a 
compatible finite element method. 
To extend the Charney-Phillips grid staggering in the compatible finite element 
framework, the buoyancy is chosen from the function space which has the same 
degrees of freedom as the vertical part of the velocity space. 
As the buoyancy space is discontinuous in the horizontal direction and continuous 
in the vertical direction, we proposed a blend of an upwind DG method in the 
horizontal direction and SUPG method in the vertical direction.

The model reproduced several quasi-periodic lifecycles of fronts despite the 
presence of strong discontinuities. 
The general results of frontogenesis are consistent with the early studies. 
To validate the numerical implementation and assess the long term performance 
of the model, the asymptotic convergence to the SG limit solution is examined. 
Despite the large difference with the SG solution in RMSV, the solutions of the 
vertical slice model were converging to a solution in geostrophic balance 
as the Rossby number was reduced. 
With off-centring, the model showed the expected second-order 
convergence rate from the early stage up to the formation of the discontinuity, 
and showed a first-order rate for some time afterwards.
This result is comparable to the previous results 
using a finite difference semi-implicit semi-Lagrangian method
\citep{visram2014framework}, 
indicating that the compatible finite element method 
is performing well for this test problem. In particular, the use of 
Eulerian advection schemes rather than semi-Lagrangian schemes does not
degrade the results.

The energy analysis showed that the model suffers from dissipation of kinetic 
energy of the cross-front velocity due to the lack of resolution 
at the fronts in the $v$ advection scheme. However, the energy loss is very 
small compared to the amplitudes of potential energy and kinetic energy of 
the cross-front velocity, and is unlikely to account for the large gap between 
the RMSV values of the model and the SG limit. The large gap corresponds more 
likely to the fact that the lack of resolution shuts off the growth of the RMSV of the model
about two days early compared to that of the SG limit.  
As the frontal discontinuity reaches the grid scale very quickly, we would need 
resolution of several orders of magnitude greater than presently used to reach 
the RMSV of the SG solution. 

The frontogenesis test case shown in this paper demonstrates several aspects 
of the compatible 
finite element framework including the treatment of advection terms in the
velocity equation, and an advection scheme for the
vertically-staggered temperature space proposed here. In concurrent
research we are incorporating these techniques into a discretisation
of the compressible Euler equations for NWP, and the formulation and
test cases will be reported in a future paper.  For a scalable
solution approach for the implicit linear system in
\eqref{eq:R_u_w}--\eqref{eq:R_p_sigma}, we are currently developing a
hybridisation capability within the Firedrake package and intend to
use this in future versions of the code. As an extension of the 
vertical slice modelling of nonlinear Eady waves, we are also considering a 
development of a parameterisation scheme which could prevent the model 
shutting off the growth of the RMSV too early, so that we could increase the 
peak amplitudes of the fronts without using unrealistically high resolution.

\section*{\large Acknowledgement}
We would like to acknowledge NERC grants NE/K012533/1 and
NE/K006789/1, EPSRC Platform grant EP/L000407/1, and the Firedrake
Project. LM additionally acknowledges support from EPSRC grant
EP/M011054/1. The source code for this project is located at
\url{https://bitbucket.org/colinjcotter/slicemodels} with the tag 
\url{paper.20161110}; 
the simulation code itself is located within this repository under
\url{paper_examples/}.  This version of the code is also archived on
Zenodo: \cite{zenodo_eady_code}.

All numerical experiments in this paper were performed with the
following versions of software, archived on Zenodo:
\cite{zenodo_firedrake};
\cite{zenodo_pyop2};
\cite{zenodo_tsfc};
\cite{zenodo_fiat};
\cite{zenodo_ufl};
\cite{zenodo_coffee};
\cite{zenodo_petsc};
\cite{zenodo_petsc4py}.
\bibliography{EadyFEM}

\end{document}